\documentclass[11pt, epsfig]{article}
\usepackage{epsfig, amsmath, amssymb, amsthm, times}
\newtheorem {thm}{Theorem}[section]
\newtheorem {lem}[thm]{Lemma}

\theoremstyle{defintion}
\newtheorem {df}[thm]{Definition}

\theoremstyle{remark}
\newtheorem{rem}[thm]{Remark}

\def\qed{\hfill $\Box$ \hfill \\}

\def\R{{\mathbb R}}

\def\lbl{\label}
\def\be{\begin{equation}}
\def\ee{\end{equation}}
\def\lbl{\label}

\def\t{\mathsf{T}}

\def\pf{\noindent{\em Proof:\ }}

\def\E{{\mathbb E}}
\def\P{{\mathbb P}}
\def\eps{\epsilon}

\def\tr{{\mathbb Tr}}

\title{The Poincare map of randomly perturbed periodic motion}
\author{
Pawel Hitczenko and
Georgi S. Medvedev\thanks{
Department of Mathematics, Drexel University, 3141 Chestnut Street,
Philadelphia, PA 19104,
{\tt \{phitczen, medvedev\} @drexel.edu} }
}

\begin{document}
\maketitle
\begin{abstract} 
A system of autonomous differential equations with a stable limit cycle
and perturbed by small white noise is analyzed in this work. In the vicinity
of the limit cycle of the unperturbed deterministic system, we define, construct,
and analyze the Poincare map of the randomly perturbed periodic motion.
We show that the time of the first exit from a small neighborhood of
the fixed point of the map, which corresponds to the unperturbed periodic
orbit, is well approximated by the geometric distribution. 
The parameter of the geometric distribution tends zero together with the 
noise intensity. Therefore, our result can be interpreted as an estimate
of stability of periodic motion to random perturbations. 

In addition, we show that the geometric distribution of the first exit 
times translates
into statistical properties of solutions of important differential equation
models in applications. To this end, we demonstrate three  examples
from mathematical neuroscience featuring complex oscillatory patterns
characterized by the geometric distribution. We show that in each of these
models the statistical properties of emerging oscillations are fully explained
by the general properties of randomly perturbed periodic motions
identified in this paper. 

\vspace{0.2cm}
\end{abstract}

\section{Introduction}
Accurate description of many important dynamical phenomena in science and engineering
is impossible without taking into account random factors. Even small random perturbations
can transform deterministic dynamics in unexpected ways and create new asymptotic
regimes, which are not present in the unperturbed deterministic system. Examples
include large deviation type mechanisms of regular dynamics in randomly perturbed systems
\cite{FR01,HM,MVE05}, stochastic resonance \cite{BSV, BG}, stochastic stabilization 
\cite{Khasminsky, Mao94}, and noise-induced synchronization \cite{GP}, to name a few.
Mathematical analysis of these and other related phenomena requires effective geometric
theory of randomly perturbed dynamical systems, which belongs to the interface
between two mathematical disciplines: the theories of dynamical systems and stochastic
processes.  The goal of this work is to extend the Poincare map method, the main geometric 
tool for studying stability of a periodic motion in deterministic systems to solutions
of randomly perturbed differential equations. 

The mathematical analysis of effects of random perturbations on dynamics of nonlinear
systems was initiated by Pontryagin, Andronov, and Vitt in their pioneering  
paper \cite{PAV}. The first systematic investigation of stability of solutions
of stochastic differential equations was undertaken by Khasminsky, who extended
many methods of classical theory of ordinary differential equations 
(cf. \cite{Malkin, hale_odes}) to randomly 
perturbed systems \cite{Khasminsky}. Freidlin and Wentzell developed the asymptotic 
method of analysis of randomly perturbed dynamical systems based on large deviations
estimates \cite{FW}. Asymptotic properties of solutions of stochastic ordinary 
differential equations were studied in \cite{Skor, Fri_sodes}. More recent approaches
for studying randomly perturbed dynamical systems are based on  ideas from
the theory of dissipative dynamical systems \cite{Arn} and those from the geometric
theory for slow-fast systems \cite{BG}. A survey of asymptotic methods for randomly 
perturbed dynamical systems with a variety of applications is available in \cite{SHS}.

In qualitative theory of nonlinear dynamical systems, local stability analysis of 
invariant sets such as equilibria and periodic orbits plays an important role.
There are many effective analytical techniques for studying stability of solutions
of deterministic differential equations \cite{hale_odes}. For randomly perturbed
systems, stability of invariant sets is reflected in the statistics of the 
times of the first exit from the corresponding domains. For domains containing
a stable fixed point, the Freidlin-Wentzell theory of large deviations characterizes
the asymptotics of the first exit time and the geometric location of the point of
exit a random trajectory from the domain \cite{FW}. Furthermore, it is known that
the limiting distribution of the first exit time is exponential \cite{Day83}.  

In the hierarchy of invariant sets of autonomous differential equations, equilibria 
are followed by periodic orbits. The main tool for studying stability of a periodic orbit
is the Poincare map (PM). The PM captures the behavior of trajectories in a typically
small neighborhood of the periodic orbit. The fixed point of the Poincare map corresponds
to the periodic orbit. Stability of this fixed point of the map translates into the 
stability of the periodic orbit. In this work, we consider
an autonomous system of differential equations in $\R^n$ possessing a limit
cycle and perturbed by state-dependent white noise process. Under general assumptions
on the limit cycle, we derive the PM for the randomly perturbed system and study its 
stability. If the periodic orbit
of the deterministic system is asymptotically stable, we show that the distribution
of the first exit times is approximately geometric with the parameter of the 
geometric distribution tending to zero together with the noise intensity. 
This result can be interpreted as a form of stability of motion near the periodic orbit:
while the trajectories of the PM eventually leave a small neighborhood of the fixed point,
for small noise they remain in this neighborhood for a very long time with 
large probability. 

The geometric distribution of the first exit times resulting from stochastic perturbations of 
a stable limit cycle is important for understanding statistical properties of nonlinear 
oscillations generated by the randomly perturbed models. In fact, this work was
motivated by our earlier analysis of a class of neuronal models in \cite{HM}. In conclusion
of this paper, we discuss the applications of our results to irregular bursting and 
mixed-mode oscillations, two important oscillatory regimes encountered in conductance
based models of neurons. 

The organization of the paper is as follows.  Section~\ref{assumptions} contains the formulation
of the problem, preliminaries about the PM, and the statement of the main result. 
In Section~\ref{sec.Pmap}, in a small neighborhood of the periodic orbit of the deterministic
system, we derive a PM for the randomly perturbed problem. The analysis follows the derivation
of the PM for deterministic ordinary differential equations \cite{hale_odes} using
the asymptotic expansions of solutions of stochastic differential equations on finite
time intervals \cite{BL,BL62,FW}. 
Section~\ref{sec.linearization} presents the analysis
of the linearization of the PM. Here, we use the results of Kesten on the iterations 
of linear functions of random matrices (cf.~\cite{kes}) to show that the first exit time
of the linearized problem has asymptotically geometric distribution. The analysis of this
section is a substantial generalization of our previous work \cite{HM}. In Section~\ref{nlin},
we estimate the contribution of the nonlinear terms of the PM on the statistics of the 
first exit times. This concludes the proof our main result - Theorem~\ref{full}.
In Section~\ref{example}, we illustrate the analysis of the previous sections
with applications to three problems in mathematical neuroscience. 
To this end, we use two conductance-based models of single neurons
and a model of electrically coupled network of pancreatic $\beta-$cells.
In the absence of noise, all these models exhibit stable periodic 
oscillations. Adding noise to these models results in more complex 
stochastic oscillatory regimes: irregular bursting and mixed-mode 
oscillations. We show that despite different mathematical formulations 
of these models and different forms of resultant oscillations,
emergent stochastic oscillatory patterns in these models are formed
due to the random perturbations of stable limit cycle oscillations.
We verify numerically that the number of spikes in one burst and 
the number of small oscillations between consecutive spikes in the time series
corresponding to these models are distributed approximately geometrically
in accord with Theorem~\ref{full}. 
We conclude this paper with a brief discussion of results of
this work as well as related results in the literature in Section~\ref{discuss}.

\section{Assumptions and results}\lbl{assumptions}

\setcounter{equation}{0}
\subsection{The model} \lbl{the-model}
Consider an ordinary differential equation for $x:\R\rightarrow\R^{d+1}$
\be\lbl{ode}
\dot x=f(x),
\ee
where $f:\R^{d+1}\rightarrow\R^{d+1}$ is twice continuously differentiable function,
and all its partial derivatives up to the second order are uniformly bounded 
in $\R^{d+1}$.  Suppose 
\be\lbl{period}
x=u(t),\; u(t+1)=u(t),\; t\in\R
\ee
is a nonconstant periodic solution of (\ref{ode}) with the least period $1$
and nonvanishing derivative
\be\lbl{nonvanish}
\dot u(t)\neq 0, \; t\in\R.
\ee
The corresponding orbit is denoted by
\be\lbl{orbit}
\mathcal{O}=\{x=u(\theta),\; \theta\in \mathbb{S}^1:=\R^1/\mathbb{Z} \}.
\ee
Along with (\ref{ode}) we consider a randomly perturbed system 
\be\lbl{sode}
\dot x_t=f(x_t) +\sigma P(x_t)\dot W_t,
\ee
where $W_t$ is a standard Wiener process in $\R^{d+1}$. Matrix $P(x)\in\R^{(d+1)\times (d+1)}$
is nondegenerate for any $x\in \R^{d+1}$. The entries of $P$ are continuously differentiable 
functions and all their
partial derivatives are uniformly bounded in $\R^{d+1}$.
The noise intensity $\sigma\ge 0$ is considered a small parameter.
Equation (\ref{sode}) is understood in the sense of Ito \cite{o}.

\subsection{The local coordinates}
By Theorem VI.1.1 of \cite{hale_odes}, in a small neighborhood of $\mathcal{O}$, 
there is an orthonormal moving coordinate frame 
\be\lbl{basis}
\{ v(\theta), z_1(\theta), z_2(\theta),\dots, z_{n-1}(\theta)\},
\quad v(\theta)={\dot u(\theta)\over \left|\dot u(\theta)\right|},\;
\;\theta\in\mathbb{S}^1.
\ee
such that the first vector $v(\theta)$ points in the tangential direction to $\mathcal{O}$.
Denote
$$
Z(\theta)=\mbox{col}(z_1(\theta),\dots, z_{n-1}(\theta))\in\R^{(d+1)\times d},
$$
then
\be\lbl{localvar}
x=u(\theta)+Z(\theta)\rho,\; 
\ee
defines a smooth invertible transformation 
$x\mapsto (\theta,\rho)\in \mathbb{S}^1\times \R^d$ in a small neighborhood of 
$\mathcal{O}$ (cf.~Theorem VI.1.1 of \cite{hale_odes}).

We use the moving coordinates to rewrite (\ref{ode}) in a form more amenable for
analysis.
\begin{lem}\lbl{local-coordinates}
In new coordinates (\ref{localvar}), 
near $\mathcal{O}$ (\ref{sode})
has the following form
\begin{eqnarray}\lbl{theta}
\dot\theta_t &=& 1+a(\theta_t)^\t\rho+O(\left|\rho_t\right|^2) +
\sigma\left( h(\theta_t)^\t+O(|\rho_t|)\right)\dot W_t,\\
\lbl{eqn-rho}
\dot\rho_t &=& R(\theta_t)\rho_t +O(\left|\rho_t \right|^2)+
\sigma\left( H(\theta_t)^\t+O(|\rho_t|)\right)\dot W_t,
\end{eqnarray}
where
\begin{eqnarray}\lbl{at}
a^\t (\theta) &=& 2{v^\t(\theta)\over \left|f(u(\theta))\right|} \left(Df(u(\theta))\right)^s Z(\theta),\\
\lbl{Atheta}
R(\theta)      &=& Z^\t(\theta) Df(u(\theta))Z(\theta)-Z^\t(\theta)Z^\prime(\theta),\\
\lbl{h1}
h(\theta)&=&{1\over |f\left(u(\theta)\right)|^2} 
P\left(u(\theta)\right)^\t f\left(u(\theta)\right)\in \R^n,\\
\lbl{H1}
H(\theta) &=& P\left(u(\theta)\right)^\t Z(\theta)\in\R^{(d+1)\times d}.
\end{eqnarray}
where $M^s:=2^{-1}(M+M^\t)$ stands for symmetric part of matrix $M$.
\end{lem}

\begin{rem} By setting $\sigma=0$  on the right hand side of 
(\ref{at}) and (\ref{Atheta}), we obtain the deterministic
equation (\ref{ode}) rewritten in the moving coordinate frame
near $\mathcal{O}$.
\end{rem} 

\subsection{The Poincare map}
Our next goal is to construct the PM near the periodic orbit $\mathcal{O}$.
We first review properties of the PM of the deterministic system
(\ref{ode}) and then turn to the PM for the randomly perturbed system (\ref{sode}).

Throughout this subsection, we will use
\be\lbl{cross}
S_{0,\tilde\delta}=\left\{ x=u(0)+Z(0)\rho:\; |\rho| \le \tilde\delta \right\},
\ee
to denote a local section to periodic orbit $\mathcal{O}$. Here, $\tilde\delta>0$ is 
assumed to be sufficiently small. 

Using the continuity of solutions of (\ref{ode}) with respect to initial data,
for sufficiently small $\delta_1>0$, one can find $\delta_2\ge\delta_1$ such that
for any initial condition $x_0:=(0,\rho_0)\in S_{0,\delta_1}$, the first return time
\be\lbl{first-time}
T(x_0)=\inf \{t>0:\; x(t)\in S_{0,\delta_2}\}
\ee
is finite. The PM $\mathbf{P}:~S_{0,\delta_1}\to S_{0,\delta_2 }$
is defined as follows
\be\lbl{EqnPmap}
\mathbf{P}(\rho_0):=\rho_T,
\ee
where $\rho_{\{0,T\}}=Z(0)^{-1}x_{\{0,T\}}$. 

The following properties of $\mathbf{P}$ are well-known (see, e.g., \cite{SSTC}):
\begin{description}
\item[A)]  $\rho=0$ is a fixed point of $\mathbf{P}$. 
\item[B)]
For small $|\rho|,$ the Poincare map has the following form
\be\lbl{linear}
\mathbf{P}(\rho)=A\rho + O(|\rho|^2), \; A=X(1),
\ee 
where $X(t)$ is the principal matrix solution of the homogeneous system
\be\lbl{homo}
\rho = R(t) \rho \qquad (\mbox{cf.}~ (\ref{eqn-rho})\;\mbox{and}\;(\ref{Atheta})).
\ee
\item[C)]
If the moduli all eigenvalues of $A$, $\mu_1,\mu_2, \dots, \mu_{n-1}$, are less than $1$,
then periodic orbit $\mathcal{O}$ is asymptotically orbitally stable.
\end{description}

Next, we turn to the randomly perturbed problem
(\ref{sode}).
Let  $x(t)$ and $x_t:=(\theta_t, \rho_t)$ denote  the solutions of
the deterministic and randomly perturbed problems (\ref{ode}) and (\ref{sode}) 
respectively starting with with initial condition $x_0=(0,\rho_0)\in S_{0,\delta_1}$.
Using the large deviation estimates for solutions of (\ref{sode}) 
(cf.~Lemma~2.1 of Chapter~4 in \cite{FW}\footnotemark), we have 
\be\lbl{large-deviations}
\P_{x_0}(\sup_{t\in [0,2]} \left|x_t-x(t)\right|\ge \delta_2)\le \exp\{-c\sigma^{-2}\},
\ee
for certain constant $c>0$, provided $\sigma>0$ is sufficiently small.
\footnotetext[1]{Lemma~2.1 in \cite{FW} is stated for a system of the form (\ref{sode})
with $P(x)=\mbox{id}$. However, the argument used in the proof of this lemma applies 
to systems with $P(x)$ satisfying the assumptions in \S\ref{the-model}
after a suitable modification of the action functional.}
Therefore, with exponentially close to $1$ probability, the trajectory of (\ref{sode})
$x_t$ returns to $S_{0,2\delta_2}$ after the first return time
\be\lbl{stop}
\tau_\sigma=\inf\{t>0.5:\; \theta_t=1\}.
\ee
The first return map $\bar \rho=\mathbf{P}(\rho)$ is defined by
\be\lbl{map}
\mathbf{P}(\rho):=\rho_{\tau_\sigma}.
\ee

The following lemma provides the asymptotic description of the PM.
\begin{lem}\lbl{local-map}
The first return map has the following form
\be\lbl{rmap}
\bar\rho=A(I+\sigma B\zeta)\rho+\sigma\eta +O(\sigma^2, |\rho|^2),
\ee
where 
\be\lbl{matrices}
A=X(1),\quad B=X^{-1}(1)\dot X(1),\quad
\ee
\be\lbl{etc}
\zeta:=\xi-b(1)^\t\eta, \; \xi=\int_0^1 h(s)^\t dW_s,\; \eta=\int_0^1 X(t)X^{-1}(s)H(s)^\t dW_s,
\ee
and
\be\lbl{bs}
b(t)^\t=-\int_0^t a(s)^\t X(s)X^{-1}(t)ds.
\ee
\end{lem}
\begin{rem}
The principal matrix solution of the system with periodic coefficients (\ref{homo})
can be written as $X(t)=Q(t)\exp\{tA\}$, where $Q(t)$ is a $1-$periodic matrix. Therefore,
$B$ in (\ref{matrices}) can be rewritten as follows
\be\lbl{rewriteB}
B=Q^{-1}(0)\dot Q(0)+Q^{-1}(0)AQ(0).
\ee
\end{rem}

\subsection{Stability of the randomly perturbed PM}

If the moduli of all eigenvalues of $A$ (cf. (\ref{EqnPmap})), 
$\mu_1,\mu_2, \dots, \mu_d$, are less than $1$, then the periodic orbit $\mathcal{O}$ 
of the deterministic system (\ref{ode}) is asymptotically orbitally stable.
This means that if the initial condition $x(0)$ is chosen sufficiently close to 
$\mathcal{O}$, the trajectory $\{x(t), t\ge 0\}$ will remain the vicinity of $\mathcal{O}$
for all future times. 
Even for very small noise intensity $\sigma>0$, for any initial condition
a generic trajectory of the randomly perturbed system (\ref{sode}) eventually leaves any 
neighborhood of $\mathcal{O}$ due to the large deviations  (cf. \cite{FW}).
Nonetheless on finite interval of time $t\in [0,T]$, $x_t$ 
with high probability exhibits stable behavior,
provided $\sigma>0$ is sufficiently small and $A$ is a stable matrix.
Thus, we expect that the trajectory of the PM remains close to the 
origin for a long time. 
To describe stability properties of the randomly perturbed PM,
for the trajectory of (\ref{rmap}) $\{\rho_n\}$ starting in a measurable set $D\subset\R^d$
containing the origin, we define the first exit time
\be\lbl{exit}
\tau(\rho_n, D) =\min\{ n:\; \rho_n\notin D\}.
\ee

In the remainder of this subsection, we formulate two theorems characterizing
the distribution of the first exit time of the
trajectories of the linearized PM (Theorem~\ref{asymptotics})
and those of the full nonlinear map (Theorem~\ref{full}).
These theorems use certain auxiliary notation, which we review next. 

\begin{df}\lbl{def:asgeom}\cite{HM}
Let $Y$ be a random variable with values in the set of positive integers and let 
$0<p<1$. We say that $Y$ is asymptotically geometric with parameter $p$ if
\begin{equation}\lbl{asgeom}
\lim_{n\to\infty} \frac{\P\left(Y=n\right)}{\P(Y\ge n)}= p.
\end{equation}
\end{df}

Recall the expression of the PM in (\ref{rmap}). Along with the full nonlinear
PM (\ref{rmap}), we  also consider its linearization  
\be\lbl{linearization}
\varrho_n= A(I+\sigma B\zeta_n)\varrho_{n-1} +\sigma\eta_n,
\ee
where $\{\zeta_n\}$ and $\{\eta_n\}$ are Gaussian random processes defined in (\ref{etc}).

Suppose $A$ is a stable matrix with the spectral radius 
\be\lbl{radius}
\mathbb{\rho}(A)<1-\epsilon
\ee
for some $0<\epsilon<1$.
Then there exists matrix norm in $\R^{d\times d}$ such that 
\be\lbl{norm}
\| A\|^\prime\le 1- \epsilon,
\ee
(see, for example, \cite{hj}). We will refer to the norm
in (\ref{norm}) as the adapted norm.

Let $|\cdot|^\prime$  
denote a vector norm in 
$\R^d$ compatible with (\ref{norm})
\be\lbl{contraction}
|Ax|^\prime\le (1-\epsilon)|x|^\prime\quad\forall x\in \R^d,
\ee
and $\gamma_2\ge \gamma_1>0$ such that 
\be\lbl{equiv}
\gamma_1 |x|\le |x|^\prime \le \gamma_2 |x| \quad \forall x\in\R^d.
\ee

For fixed  $h>0$ we
define 
\be\lbl{domain}
D_h=\{x\in\R^{d}:\; |x|^\prime\le h \}.
\ee

It is instructive first to understand the statistics of the first exit times
for the linearized PM (\ref{linearization}).
\begin{thm}\lbl{asymptotics}
Suppose (\ref{radius}) holds and $\{\varrho_n\}$ denotes the trajectory
of (\ref{linearization}) starting from initial condition $\varrho_0=0$. Then
for certain $\sigma_0=O(\epsilon)$ and for all $0<\sigma<\sigma_0$, 
$\tau (\varrho_n,D_h)$ is an asymptotically 
geometric RV with parameter
\be\lbl{parameter}
C_1\sigma \exp\left\{ {-C_2\over\sigma ^2}\right\} \left(1+O(\sigma^2)\right)
\le p \le
C_3\sigma \exp\left\{ {-C_4\over\sigma ^2}\right\} \left(1+O(\sigma^2)\right)
\ee
for some positive constants  $C_{1,2,3,4}$ independent from $\sigma$.
\end{thm}

We are now in a position to formulate the main result of this paper.
\begin{thm}\lbl{full}
Suppose (\ref{radius}) holds and $\{\rho_n\}$ denotes the trajectory
of (\ref{rmap}) starting from initial condition $\rho_0=0$.
Then for certain $\sigma_0=O(\epsilon)$ and for all $0<\sigma<\sigma_0$, 
RV $\tau (\varrho_n,D_h)$ is subject to the 
following estimate:
for $n\ge 1$,  
\be\lbl{sub}
\frac{\P(\bar\tau=n+1)}{\P(\bar\tau\ge n+1)}\le p,
\quad\mbox{where}\quad  0\le p\le C(\sigma h^2)^{2/3},
\ee
where $C$ is a positive constant independent of $\sigma$ and $h$.
\end{thm}

\section{The derivation of the Poincare map}\lbl{sec.Pmap}
\setcounter{equation}{0}
In this section, we compute the linear part of the Poincare 
map of the periodic orbit $\mathcal{O}$. 

\subsection{The variational equation}
The first step in the derivation of the PM is the asymptotic approximation 
of the solution of (\ref{theta}) and
(\ref{rho}) subject to initial condition
\be\lbl{id}
\theta_0=0\quad\mbox{and}\quad |\rho_0|<\delta.
\ee

\begin{lem}\lbl{expand}
On any finite interval of time, for sufficiently small $\sigma>0$, the solution
of the initial value problem (\ref{theta}), (\ref{rho}), and (\ref{id})
admits the following asymptotic expansion
\be\lbl{sumup}
\theta_t=t-b(t)^\t X(t)\rho_0+
\sigma\left(\zeta_t+O(|\rho_0|)\right) +O(\sigma^2, |\rho_0|^2),\;
\mbox{and}\;
\rho_t=X(t)\rho_0+\sigma\eta_t +O(\sigma^2, |\rho_0|^2),
\ee
where 
\be\lbl{theta1}
\eta_t=\int_0^t X(t,s)H(s)^\t dW_s,
\; \xi_t=\int_0^t h(s)^\t dW_s,\;
\;\;\mbox{and}\;\; \zeta_t=\xi_t-b(t)^\t\eta_t,
\ee
and 
\be\lbl{b(t)}
b(\theta)^\t=-\int_0^\theta a^\t(s) X(s)X^{-1}(\theta) ds.
\ee
\end{lem}
\pf
For integration of (\ref{theta}) and (\ref{rho}), it is convenient
to change $\theta$ to a new variable
\be\lbl{phi}
\phi=\theta +b(\theta)^\t\rho,
\ee
where $b(t)$ is defined in (\ref{b(t)}).
In new coordinates (\ref{at}) and (\ref{Atheta})
become
\begin{eqnarray}\lbl{theta}
\dot\phi_t &=& 1+O(\left|\rho_t\right|^2) +
\sigma\left( h(\phi_t)^\t+O(|\rho_t|)\right)\dot W_t,\\
\lbl{rho}
\dot\rho_t &=& R(\phi_t)\rho_t +O(\left|\rho_t \right|^2)+
\sigma\left( H(\phi_t)^\t+O(|\rho_t|)\right)\dot W_t.
\end{eqnarray}
The corresponding initial condition is
\be\lbl{id-1}
\phi_0=0\quad\mbox{and}\quad |\rho_0|<\delta.
\ee
On a finite time interval $t\in [0,2]$, we expand $x_t=(\phi_t, \rho_t)$
in the asymptotic sum
\be\lbl{Anz}
z_t=z^{(0)}_t+\sigma z_t^{(1)}+\mathcal{R}(t,\sigma),
\ee
where $z_t^{(0,1)}=(\phi^{(0,1)}_t, \rho^{(0,1)}_t).$ Function $z_t^{(0)}$
is the deterministic. The first order correction $z_t^{(1)}$ is a Gaussian
process. Below we state  the corresponding initial value problems
for $z_t^{(0,1)}$. The remainder
satisfies the following estimate (cf. Theorem 2.2, \cite{FW})
\be\lbl{remainder}
\E\{\sup_{t\in [0,2]} \left|\mathcal{R}(t,\sigma)\right|^2\} \le C \sigma^2
\ee
for some positive constant $C$. We denote the asymptotic relation in (\ref{remainder})
by $\mathcal{R}(t,\sigma)=O(\sigma^2).$

The zeroth order problem is given by
\be\lbl{zero}
\dot\phi^{(0)}_t=1+O(|\rho|^2) \quad\mbox{and}\quad 
\dot\rho^{(0)}_t=R(\phi)\rho^{(0)}_t+O(|\rho_t|^2),
\ee
subject to $x_t^{(0)}=(0,\rho_0)^\t$. By integrating (\ref{zero}), we have
\be\lbl{lead}
\phi^{(0)}_t=t+O(|\rho_0|^2), \quad\mbox{and}\quad
\rho^{(0)}_t=X(t)\rho_0+O(|\rho_0|^2).
\ee
The first order problem has the following form
 \be\lbl{first}
\dot\phi^{(1)}_t=h(t)^\t\dot W_t \quad\mbox{and}\quad 
\dot\rho^{(1)}_t=R(\phi)\rho^{(1)}_t+H(t)^\t\dot W_t,
\ee
with initial condition $x_t^{(1)}=(0,0)^\t$.
From (\ref{first}), we find
\be\lbl{next}
\xi_t:=\phi^{(1)}_t=\int_0^t h(s)^\t dW_s\quad\mbox{and}\quad
\eta_t:=\rho^{(1)}_t=\int_0^t X(t,s)H(s)^\t dW_s.
\ee
Combining (\ref{Anz}), (\ref{lead}) and (\ref{next}) and switching back to the
original variables $(\theta_t,\rho_t)$ (cf.~(\ref{phi})), we obtain (\ref{sumup}).\\
\qed

\subsection{Proof of Lemma~\ref{local-map}}
First, we estimate the time
of the first return.
\begin{lem}\lbl{time}
The first return time is given by
\be\lbl{expand-tau}
\tau_\sigma=1+b(1)^\t X(1)\rho_0-\sigma\zeta_1+o(\sigma)+O(|\rho_0|^2),
\ee
\end{lem}
\pf\;
Using the definition of the first return time  (\ref{stop}) and
the asymptotic expansion of $\theta_t$ (\ref{sumup}), we have
\be\lbl{setup}
\tau_\sigma-a_1(\tau_\sigma)\rho_0+\zeta_{\tau_\sigma}+O(2)=1, \quad \mbox{a.s.},
\ee
where $a_1(t):=b(t)^\t X(t)$ and $O(2):=O(\sigma^2, |\rho_0|^2)$.
It follows from (\ref{sumup}) that $\tau_\sigma\rightarrow\tau_0$ a.s., as $\sigma\to 0$,
where $\tau_0$ is the first return time for the unperturbed deterministic trajectory
($\sigma=0$). Thus, we write 
\be\lbl{split}
\tau_\sigma = \tau_0+\tau_1(\sigma), 
\ee
where $\tau_1(\sigma)=o(1)$. 
By plugging (\ref{split}) in (\ref{expand-tau}), we have
\be\lbl{plugin}
\tau_0+\tau_1(\sigma)-a_1(\tau_0)\rho_0 - a_1^\prime (t^\prime)\tau_1(\sigma)\rho_0
+\sigma \zeta_{\tau_0} +\sigma (\zeta_{\tau_\sigma}-\zeta_{\tau_0})+O(2)=1,
\ee
where $t^\prime$ is a real number lying between $\tau_0$ and $\tau_\sigma$.
By setting $\sigma=0$ in (\ref{plugin}), we obtain
$$
\tau_0-a_1(\tau_0)\rho_0+ O(|\rho_0|^2)=1,
$$
and, thus,
\be\lbl{tau-0}
\tau_0=1+a_1(1)\rho_0+ O(|\rho_0|^2).
\ee
Next,
$$
\tau_1(\sigma) (1-a_1(t^\prime)\rho_0)+\sigma\zeta_1 +
\sigma\{ (\zeta_{\tau_\sigma}-\zeta_{\tau_0}) + (\zeta_{\tau_0}-\zeta_1)\}=0.
$$
By the It\^{o} isometry,  two terms in the curly brackets are
$o(1)$ and $O(|\rho_0|)$ respectively.
Thus,
\be\lbl{tau-1}
\tau_1(\sigma)=-\sigma\zeta_1 + O(\sigma|\rho_0|) + o(\sigma).
\ee
The combination of (\ref{split}), (\ref{tau-0}), and (\ref{tau-1}) proves the lemma.\\
\qed

We are now in a position to compute the first return map. 
From (\ref{sumup}) and (\ref{expand-tau}), we have
\be\lbl{barrho}
\bar\rho=X(\tau_\sigma)\rho_0+\sigma\eta_{\tau_\sigma}+O(2)=
X(1+\sigma\zeta_1)\rho_0+\sigma\eta_1+O(2).
\ee
Further,
$$
X(1+\sigma\zeta_1)=X(1)(I+\sigma X^{-1}(1)\dot X(1)\zeta_1)+O(\sigma^2). 
$$
This proves (\ref{rmap}).

\section{The randomly perturbed linear map}\lbl{sec.linearization}
\setcounter{equation}{0}

In this section, we study stability of the linearized PM (\ref{linearization}). 
Theorem~\ref{asymptotics} will follow from the analysis of the following
slightly more general class of equations
\be\lbl{model}
y_{n}=M_ny_{n-1} +z_n,\quad n\ge1,
\ee
where
\be\lbl{data}
M_n=A(I+\sigma \xi_n 
B)\quad\mbox{and}\quad z_n=\delta G\eta_n.
\ee
Here, $A,B,G \in R^{d\times d}$; $\xi_n$ and $\eta_n$ are 
independent identically distributed (IID) copies of standard normal RV $\xi$
and $\eta$ in $\R$ and $\R^d$ respectively. Furthermore, 
we assume that $A$ is a stable matrix with spectral radius
subject to (\ref{radius})
and $G$ is nondegenerate. In fact, without loss of generality, we assume that
$G$ is positive definite, since otherwise, one can take $G:=(GG^\t)^{1/2}$
without changing the distribution of $z_n$.

\begin{thm}\lbl{perturbed-map} 
Let $h>0$ be fixed. For $y_0\in D_h$, denote the first exit
time of the trajectory of (\ref{model}) from $D_h$
(cf. (\ref{domain}))  by $\tau:=\tau(y_n,D_h)$. 
Then there exist positive $\delta_0=O(\epsilon)$ and  
$\sigma_0=O(\epsilon)$ such that for any $y_0\in D_h$, 
$\sigma\in (0,\sigma_0)$, and $\delta\in (0,\delta_0),$ 
the first exit time $\tau:=\tau(y_n, D_h)$ is asymptotically 
geometric RV with parameter $p$ satisfying
\be\lbl{main-estimate}
C_1\sigma \exp\left\{ {-C_2\over\sigma ^2}\right\} \left(1+O(\sigma^2)\right)
\le p \le
C_3\sigma \exp\left\{ {-C_4\over\sigma ^2}\right\} \left(1+O(\sigma^2)\right)
\ee
for some positive constants  $C_{1,2,3,4}$, which do not depend on $\sigma$ and $\delta$.
\end{thm}
\begin{rem}
Note that setting $\delta:=\sigma$ and with appropriate choice of $G$,
(\ref{model}) coincides with the expresion of the linearized PM
(\ref{linearization}). Therefore, Theorem~\ref{asymptotics} is a
particular case of Theorem~\ref{perturbed-map}.
\end{rem}

The proof of Theorem~\ref{perturbed-map} relies on the convergence results
for iterative processes due to Kesten \cite{kes}. 
Specifically, our assumptions 
on $\{M_n\}$ and $\{z_n\}$ guarantee that $\{y_n\}$ converge in distribution
to a random vector $y\in\R^d$ as $n\to\infty$.
Furthermore, $\xi$ and $z$ are independent of $y$. We verify the 
conditions of the Kesten's theorem implying the convergence 
of $\{y_n\}$ in Lemma~\ref{verify}.

\begin{lem}\lbl{verify}
Under our assumptions on the coefficients in (\ref{model}), stochastic process
$\{y_n\}$ converges in distribution to a random vector $y\in \R^d$ as $n\to\infty$,
such that $(\xi, z)$ are independent of $y$. 
\end{lem}

The gist of the proof of Theorem~\ref{perturbed-map} lies in estimating 
$\P(My+z\notin D_h)$. This is the subject of the following lemma. 

\begin{lem}\lbl{gist} There exist constants $0<C_5<C_6<1$
such that, 
uniformly over $s\in D_h$
\be\lbl{estimateP}
C_5 \le\P(Ms
+z\notin D_h)\le C_6.
\ee
\end{lem}

We first prove Theorem~\ref{perturbed-map} using Lemma~\ref{gist}
and then  give
the more technical proof of Lemma~\ref{gist} followed by the 
proof of auxiliary Lemma~\ref{verify}.

\pf (Theorem~\ref{perturbed-map})
We need to show
\be\lbl{ratio}
\lim_{n\to\infty} {\P(\tau=n)\over \P(\tau\ge n)}=p>0
\ee
and estimate $p$ in terms of the coefficients of Equation (\ref{model}).

First, we note
\begin{eqnarray}\nonumber
\P(\tau=n+1) &=& \P(y_{n+1}\notin D_h, y_k\in D_h, k\in [n])\\
\nonumber
&=&
\P(y_{n+1}\notin D_h| y_k\in D_h, k\in[n]) \P(y_k\in D_h, k\in[n])\\
\lbl{first-stub}
&=&
\P(y_{n+1}\notin D_h| y_k\in D_h, k\in[n])\P(\tau\ge n+1).
\end{eqnarray}
Further, using (\ref{model}) and (\ref{data}) we have 
\be\lbl{Markov}
\P(y_{n+1}\notin D_h| y_k\in D_h, k\in[n])=
\P(M_{n+1}y_n+z_n\notin D_h| y_n\in D_h).
\ee
The combination of (\ref{first-stub}) and (\ref{Markov}) shows that
(\ref{ratio}) is equivalent to convergence of
\be\lbl{ratio2}
\P(M_{n+1}y_n+z_n\notin D_h| y_n\in D_h)={\P(M_{n+1}y_n+z_n\notin D_h, y_n\in D_h)\over
\P(y_n\in D_h)} 
\ee
to a nonzero limit.

By Lemma~\ref{verify},
$$
(y_n,\xi_n, z_n)\stackrel{d}{\rightarrow} (y,\xi, z),\; n\to\infty,
$$ 
where $(\xi, z)$ is independent of $y$. 
Thus, 
the limit
of the numerator in (\ref{ratio2}) as $n\to\infty$ is 
\be\lbl{numerator}
\P(My +z\notin D_h, y\in D_h)=
\int_{D_h} \P(Ms+z\notin D_h) dF_y(s),
\ee
where $F_y(\cdot)$ stands for the distribution function of $y$.
By Lemma~\ref{gist},
$$ 
C_5 \P(y\in D_h)\le \int_{D_h} \P(Ms+z\notin D_h) dF_y(s)\le C_6
\P(y\in D_h), $$ 
and, therefore, 
the limit, $p:=\P(My+z\notin D_h|y\in D_h)$,  on the right--hand side of \eqref{ratio2}  
satisfies
$$
0<C_5\le p
\le C_6<1.
$$
\qed

\pf (Lemma~\ref{gist})
\begin{enumerate}
\item
Let $s\in D_h$, i.e., $|s|^\prime\le h$. Then, 
 using (\ref{contraction}) and (\ref{data}), we have
\be\lbl{1st-step}
\P(|Ms+z|^\prime\ge h)
\le \P(\sigma |\xi| |ABs|^\prime+\delta |G\eta|^\prime
> \epsilon h). 
\ee
Recall that $\xi\in \mathcal{N}(0,1)$ and define event
\be\lbl{defineF}
\mathcal{F}=\left\{|\xi|\le 
{\epsilon\over 2\sigma \|AB\|^\prime}\right\}.
\ee
Using (\ref{defineF}), we take the estimate in (\ref{1st-step}) one step 
further
\be\lbl{2nd-step}
\P(|Ms+z|^\prime\ge h) 
\le 
\P\left(|G\eta|^\prime>{\epsilon h\over 2\delta}, \mathcal{F}\right) +\P(\mathcal{F}^c).
\ee
We bound the second term on the right hand side of (\ref{2nd-step}),
using the normal distribution of $\xi\in\mathcal{N}(0,1)$
\begin{eqnarray}\nonumber
\P(\mathcal{F}^c) &=&\P\left( |\xi| >{\epsilon\over 2\sigma \|AB\|^\prime}\right)
\le \P\left( |\xi|^\prime >{\epsilon\over 2\sigma
    \|AB\|^\prime}\right)\\
\lbl{2nd-term}
&\le &\sqrt{{2\over\pi}}{\sigma\over C_7\epsilon}
\exp\left\{ {-C_7^2\epsilon^2\over 2\sigma^2}\right\}
\left( 1+ O(\sigma^2\epsilon^{-1})\right),
\end{eqnarray}
where $C_7=(2\| AB\|^\prime)^{-1}$ is independent from $\sigma$.
\item
Our next goal is to bound the 
first term on the right hand side
of (\ref{2nd-step}).
Gaussian random vector 
\be\lbl{z-tilde}
\tilde z= G\eta
\ee
has zero mean and 
covariance matrix $GG^\t$. 
Denote the median of the distribution
of $\tilde z$ by $m$ and
\be\lbl{concentrate}
K=\{x\in\R^d:~ |Gx|\le m
\}\quad \mbox{and}\quad
K_r=\{ x\in\R^d:~ \mbox{dist}~(x,K)\le r\}.
\ee
We will need the following concentration inequality for 
Gaussian random vectors (cf.~
\cite[pp. 19--21]{lt})
\be\lbl{concentrate}
\P(z\in K_r^c)\le 1- \Phi(r)={ \exp\{-2^{-1}r^2\}\over 
\sqrt{2\pi} r}\left(1+O(r^{-2})\right),
\ee
where $\Phi(r)=\frac{1}{\sqrt{2\pi}}\int^r_{-\infty} e^{-s^2/2}ds$. 
We also employ an observation of Kwapie\'n \cite{k}
to bound the median of $\tilde z$
\be\lbl{Kwapien}
m^2\le \E~|G\eta|^2=\E~\tr\{(G\eta)^\t(G\eta)\}= \tr~\E~(G\eta\eta^\t G^\t)
=\|G\|^2_F,
\ee 
where $\|G \|_F=\sqrt{ \tr~(GG^\t)}$ is the Frobenius norm of $G$.
From 
(\ref{Kwapien}) and the triangle inequality, 
for any $y\in K_r$ 
there exists $x\in K$ such that
\be\lbl{triangle}
|Gy|\le |Gx|+|G(y-x)|\le \|G\|_F+\|G\|r.
\ee

We are now in a position to use the concentration inequality (\ref{concentrate})
to bound the first term on the right hand side of (\ref{2nd-step}).

First, using (\ref{equiv}) we switch to the Euclidean vector norm 
\be\lbl{switch}
\P\left( |\tilde z|^\prime > {\epsilon h\over 2\delta}\right)\le
\P\left( |\tilde z| > {\epsilon h\over 2\gamma_2\delta}\right).
\ee
Next, we specify the upper bound on $\delta$
\be\lbl{specify-delta}
0<\delta<{\epsilon h\over 4 \|G\|_F \gamma_2}
\ee 
where $\gamma_2$ is given by \eqref{equiv}, and choose 
\be\lbl{pick-r}
r={\epsilon h\over 4\delta \gamma_2 \|G \|}.
\ee
Then, \eqref{specify-delta} implies that 
\[\|G\|_F\le\frac{\epsilon h}{4\delta\gamma_2},\]
and \eqref{pick-r} is
\[\frac{\epsilon h}{4\delta\gamma_2}=\|G\| r.\]
Hence,
\be\lbl{re-express}
\frac{\epsilon h}{2\delta\gamma_2}=\frac{\epsilon h}{4\delta\gamma_2}+\frac{\epsilon h}{4\delta\gamma_2}
\ge\|G\|_F+\|G\| r.\ee
Thus, using 
(\ref{concentrate})  and (\ref{re-express}), 
we obtain that the right--hand side of \eqref{switch} is bounded by
\be\lbl{one-more}
\P\left( |\tilde z| > {\epsilon h\over 2\gamma_2\delta}\right)\le 
\P\left( |\tilde z| > \|G\|_F+\|G\| r\right)\le \P(K_r^c)\le {\exp\{-2^{-1}r^2\}\over \sqrt{2\pi} r } 
\left(1+O(r^{-2})\right).
\ee
Finally, by combining (\ref{switch}) and (\ref{one-more}), recalling the definition
of $\tilde z$ in (\ref{z-tilde}), 
we arrive at
\be\lbl{almost-there}
\P\left( |G\eta|^\prime > {\epsilon h\over 2\delta}\right)
\le {\exp\{-2^{-1}r^2\}\over \sqrt{2\pi} r } \left(1+O(r^{-2})\right),
\ee
where $r$ is specified in (\ref{pick-r}).
\item
Combining (\ref{2nd-step}), (\ref{2nd-term}), and (\ref{almost-there}),
we have 
\begin{eqnarray*}
\P(M
s+z\notin D_h)&\le& 
\sqrt{{2\over\pi}}{\sigma\over C_7\epsilon}
\exp\left\{ {-C_7^2\epsilon^2\over 2\sigma^2}\right\} 
\left( 1+ O(\sigma^2\epsilon^{-1})\right)\\
&+&
{\delta\over \sqrt{2\pi} C_8\epsilon h}
\exp\left\{ -{C_8^2\epsilon^2 h^2\over 2\delta^2}\right\}\left(1+
O(\delta^2\epsilon^{-2} h^{-2})\right)
\end{eqnarray*}
\item
To get a lower bound on $\P(M
s+z\notin D_h)$, we invoke the 
Anderson inequality \cite{And}
\be\lbl{Anderson}
\P\left( (z-x)\in D_h\right)\le \P(z\in D_h) \;\;\forall x\in \R^d,
\ee
which is applicable since $D_h$ is convex and symmetric 
about $0$. By (\ref{Anderson}),
\be\lbl{by-andreson}
\P(M
s+z\notin D_h)\ge \P(|G\eta|\ge h\delta^{-1}).
\ee
Further,
\begin{eqnarray}\nonumber
\P(|G\eta|\ge h\delta^{-1}) &\ge& \P(|\eta|\ge h(\lambda_1(G)\delta)^{-1})\ge \\
\lbl{further}
\P\left(|z|\ge h\left(\lambda_1(G)\delta\right)^{-1}\right)&=&\sqrt{{2\over \pi}} 
{\delta\over C_9 h}
\exp\left\{ {-C_9^2 h^2\over 2\delta ^2}\right\} \left( 1+O(\delta^2h^{-2})\right).
\end{eqnarray} 
\end{enumerate} 
\qed

\noindent \pf~(Lemma~\ref{verify})
Recall that according to Kesten's results (see the beginning of Section~3
or Theorem~6 in \cite{kes}) for the convergence in distribution of
$\{y_n\}$ we need to know that  
\[\E|G\eta|^\gamma<\infty\mbox{\ for some $\gamma>0$, 
 and that\ }
\alpha:=\lim_{n\to\infty}\|\prod_{k=1}^nA(I+\sigma\xi_kB)\|<0.\]
(The existence of such $\alpha$ is guaranteed by the condition
$\E\log^+\|A(I+\sigma\xi B)\|<\infty$, see \cite{fk}.)
The first of these conditions is obviously true and holds for every
$\gamma>0$. To verify the  second condition, we next show 
\be\lbl{second}
\alpha:=\lim_{n\to\infty}\frac1n\log\|\prod_{k=1}^nA(I+\sigma\xi_kB)\|<0,
\ee
provided 
\be\lbl{sigma}
\sigma<\frac{\epsilon}{c(1-\epsilon)}.
\ee
for a constant $c>0$ specified below. 

 As we mentioned, the existence of $\alpha$ is guaranteed by the
integrability of $\log\|A(I+\xi B)\|$. Furthermore, by
\cite[Theorems~1 and~2]{fk} it suffices to show that 
\be\lbl{lln}
\lim_{n\to\infty}\frac1n\E\log\|\prod_{k=1}^nA(I+\sigma\xi_kB)\|<0.
\ee
Obviously, we can replace $\|\ \cdot\ \|$ by any other 
matrix norm 
since neither condition is affected by passing to an equivalent norm.

We will show that \eqref{lln} holds for 
 the norm $\|\cdot\|^\prime$.
Since $x\to\log x$ is concave, by Jensen's inequality
\[\E\log\|\prod_{k=1}^nA(I+\sigma\xi_kB)\|^\prime\le
  \log\E\|\prod_{k=1}^nA(I+\sigma\xi_kB)\|^\prime.\]
Using corresponding of 
the adapted norm, independence of $\xi_k$'s and
triangle inequality we see  that 
\[\E\|\prod_{k=1}^nA(I+\sigma\xi_kB)\|^\prime\le\prod_{k=1}^n\E\|A(I+\sigma\xi_kB)\|^\prime\le
(\|A\|^\prime(1+\sigma c))^n,\]
where $c=\|B\|^\prime\E|\xi|=\|B\|^\prime\sqrt{2/\pi}$.

If $\sigma>0$ satisfies \eqref{sigma} then  \eqref{norm} implies
that $\|A\|^\prime(1+\sigma c)<1$ for some $c>0$. Thus, by \eqref{lln} the
\eqref{second} holds with $\alpha\le\log((1-\eps_1)((1+\sigma c))<0$.
\qed 

\section{The nonlinear estimates}\lbl{nlin}
\setcounter{equation}{0}
In this section, we prove Theorem~\ref{full}. To this end, we recall
certain notation used in the previous sections. 
The full and linearized PMs (cf.~(\ref{rmap}) and (\ref{linearization}),
respectively) are given by the following difference equations:
\begin{eqnarray}\lbl{recall-rmap}
\rho_n&=&A(I+\sigma B\zeta_{n})\rho_{n-1}+\sigma\eta_n +r(\rho_{n-1}),\\
\lbl{recall-linearization}
\varrho_n&=&A(I+\sigma B\zeta_{n})\varrho_{n-1}+\sigma\eta_n,
\end{eqnarray}
where $r(\rho)=O(|\rho|^2,\sigma^2)$.

For fixed $h>0$, we define
\be\lbl{recall-ellipsoid}
D_h=\{x\in\R^{d}:\; |x|^\prime\le h \}.
\ee
and the fist exit time $\bar \tau:=\tau(\rho_n, D_h)$.

Our goal is to prove the statement of the Theorem~\ref{full}:
for $n\ge 1$,  
$$
\frac{\P(\bar\tau=n+1)}{\P(\bar\tau\ge n+1)}\le p,
\quad\mbox{where}\quad  0\le p\le C(\sigma h^2)^{2/3}.
$$

As in the proof of Theorem~\ref{perturbed-map},
we obtain
\be\lbl{tau_bar}
\P(\bar\tau=n+1)=\frac{\P(|\bar\rho_{n+1}|^\prime>h,|\bar
 y_n|^\prime\le h)}{\P(|\bar\rho_n|^\prime\le h)}\ \P(\bar\tau\ge n+1).
\ee

Denote $x_n=\rho_n-\varrho_n$ and consider the difference equation
\be\lbl{pre-x_n}
x_n=A(I+\sigma\xi_nB)x_{n-1}+r(\rho_{n-1}).
\ee
Iterating (\ref{pre-x_n}),
we obtain
\be\lbl{x_n}
x_n=\left(\prod_{k=0}^{n-1}A(I+\sigma\xi_{n-k}B)\right)x_0+
\sum_{j=2}^{n+1}\left(\prod_{k=0}^{n-j}A(I+\sigma\xi_{n-k}B)\right)r(\rho_{j-2}),
\ee
where we used the convention on the right--hand side that the product
over the empty range is 1.   Note that  $n\le\bar\tau$ implies that
for $j=2,\dots,n+1$ we have $|\rho_{j-2}|^\prime\le h$. Hence,
there exists a constant $C$ such that 
for all such $j$'s $|r(\bar
y_{j-2})|^\prime\le Ch^2$.
Also, since  $x_0=0$ the first term
in \eqref{x_n} vanishes.  Furthermore, for $t\in\Bbb R$ and $w\in \Bbb R^d$
\[|A(I+tB)w|^\prime\le|Aw|^\prime+|A(tBw)|^\prime\le
(1-\epsilon)(|w|^\prime+|tBw|^\prime)\le (1-\epsilon)(1+|t|\cdot\|B\|^\prime)|w|^\prime,
\]
which implies that 
\[
\left|\prod_{k=0}^{n-j}A(I+\sigma\xi_{n-k}B)r(\rho_{j-2})\right|^\prime\le
(1-\epsilon)^{n-j+1}\prod_{k=0}^{n-j}(1+\sigma|\xi_{n-k}|\cdot\|B\|^\prime)|r(\rho_{j-2})|^\prime.
\]
Therefore, whenever $n\le\bar\tau$, we have
\[|x_n|^\prime\le Ch^2\sum_{j=2}^{n+1}(1-\epsilon)^{n-j+1}\prod_{k=0}^{n-j}
(1+\sigma|\xi_{n-k}|\cdot\|B\|^\prime)\le
Ch^2\sum_{\ell=0}^\infty
(1-\eps_1)^\ell \prod_{k=1}^{\ell}(1+\sigma|\xi_{k}|\cdot\|B\|^\prime).\]
where  $C=C(\eps,\sigma, A,B)$ is an absolute constant whose value may vary from use
to use.  We  bound  the variance of the latter sum as follows:
\begin{eqnarray*}
&&\operatorname{var}\left(\sum_{\ell=0}^\infty(1-\eps)^\ell
  \prod_{k=1}^{\ell}(1+\sigma|\xi_{k}|\cdot\|B\|^\prime)\right)=
\sum_{\ell=1}^\infty(1-\eps)^{2\ell}\operatorname{var}\left(\prod_{k=1}^{\ell}(1+\sigma|\xi_{k}|\cdot\|B\|^\prime)\right)\\
\quad&&+
2 \sum_{\ell=1}^\infty\sum_{m=\ell+1}^\infty(1-\eps)^{\ell+m}\operatorname{cov}\left(\prod_{k=1}^{\ell}(1+\sigma|\xi_{k}|\cdot\|B\|^\prime),\prod_{k=1}^{m}(1+\sigma|\xi_{k}|\cdot\|B\|^\prime)\right)
\\
\quad&&\le
2 \sum_{\ell=1}^\infty\sum_{m=\ell}^\infty(1-\eps)^{\ell+m}\operatorname{cov}\left(\prod_{k=1}^{\ell}(1+\sigma|\xi_{k}|\cdot\|B\|^\prime),\prod_{k=1}^{m}(1+\sigma|\xi_{k}|\cdot\|B\|^\prime)\right).\end{eqnarray*}
By independence, for $m\ge\ell\ge0$ we have
\be\lbl{prod_bound}\E
\prod_{k=1}^{\ell}(1+\sigma|\xi_{k}|\cdot\|B\|^\prime)\prod_{k=1}^{m}(1+\sigma|\xi_{k}|\cdot\|B\|^\prime)
=(\E(1+\sigma|\xi|\cdot\|B\|^\prime)^2)^\ell(1+\sigma\mu_B)^{m-\ell}
,\ee
where we have set 
\be\lbl{mub}\mu_B=\|B\|^\prime\cdot\E|\xi|=\|B\|^\prime\sqrt{2/\pi}.\ee
Therefore, the covariance above is equal to 
\begin{eqnarray*}&&(\E(1+\sigma|\xi|\cdot\|B\|^\prime)^2)^\ell(1+\sigma\mu_B)^{m-\ell}-(1+\sigma\mu_B)^{m+\ell}
\\\quad&&=(1+\sigma\mu_B)^{m-\ell}\left\{(\E(1+\sigma|\xi|\cdot\|B\|^\prime)^2)^\ell-((1+\sigma\mu_B)^2)^{\ell}\right\}.
\end{eqnarray*}
 The term in the curly brackets is bounded above by 
\be\lbl{sum_bound}(\E(1+\sigma|\xi|\cdot\|B\|^\prime)^2-(1+\sigma\mu_B)^2)\sum_{j=0}^{\ell-1}
(\E(1+\sigma|\xi|\cdot\|B\|^\prime)^2)^j((1+\sigma\mu_B)^2)^{\ell-1-j}.
\ee
The  factor outside the summation is equal to
$\operatorname{var}(\sigma\|B\|^\prime\cdot|\xi|)=O(\sigma^2)$.  To
bound the sum in \eqref{sum_bound} we
 use the inequalities
\[(1+\sigma\mu_B)^2\le
\E(1+\sigma|\xi|\cdot\|B\|^\prime)^2=1+2\sigma\|B\|^\prime\sqrt{2/\pi}+(\sigma\|B\|^\prime)^2\le (1+\sigma\|B\|^\prime)^2\]
which imply  that the sum above is bounded by 
\[\sum_{j=0}^{\ell-1}(1+\sigma\|B\|^\prime)^{2j}(1+\sigma\|B\|^\prime)^{2(\ell-1-j)}\le\ell(1+\sigma\|B\|^\prime)^{2(\ell-1)}.\]
Combining \eqref{prod_bound} and \eqref{sum_bound} with the expression
bounding  the variance we get
\be\lbl{fin_bound}
C\sigma^2\sum_{\ell=1}^\infty\ell(1-\eps)^{2l}(1+\sigma\|B\|^\prime)^{2(\ell-1)}\sum_{m=\ell}^\infty((1-\eps)(1+\sigma\mu_B))^{m-\ell}.
\ee
Using \eqref{mub} we see that if $\sigma\le\sigma_0$ where, say, $\sigma_0\le
\eps/2\|B\|^\prime$ then the double sum in \eqref{fin_bound} is bounded by a constant
depending on $\eps$  only, and
thus 
\[\operatorname{var}\left(\sum_{\ell=0}^\infty(1-\eps)^\ell
  \prod_{k=1}^{\ell}(1+\sigma|\xi_{k}|\cdot\|B\|^\prime)\right)\le
C\sigma^2.\] 
Therefore, by Chebyshev's inequality, for any $h_1>0$
we have
\be\lbl{rhs}
\P(|x_n|^\prime>h_1)\le
\P(\sum_{j=2}^{n+1}(1-\epsilon)^{n-j+1}\prod_{k=0}^{n-j}(1+\sigma|\xi_{n-k}|\cdot\|B\|^\prime)\ge\frac{h_1}{Ch^2})
\le
C\sigma^2\frac{h^4}{h_1^2}.
\ee
Note that if $h_1=O(\sigma^\alpha h^\beta)$ for $0<\alpha<1$ and
$1<\beta<3/2$. (The optimal values of $\alpha$ and $\beta$ will
be identified below.) 
Then the right hand side of (\ref{rhs}) is $O(\sigma^{2(1-\alpha)}h^{2(2-\beta)})$.
For such $h_1$, set $h_0=h-h_1$. Since $\rho_n=\varrho_n+x_n$ we have
\begin{eqnarray*}\P(|\rho_n|^\prime\le  h)&=&\P(| y_n+x_n|^\prime\le  h)\ge
\P(| y_n|^\prime\le h-h_1,|x_n|^\prime\le h_1)\\ &=&\P(|y_n|^\prime\le
h_0)-\P(|y_n|^\prime\le h_0,|x_n|^\prime> h_1)\\&\ge&\P(|y_n|^\prime\le h_0)-\P(|x_n|^\prime> h_1)
\end{eqnarray*} 
This bounds from below the denominator in \eqref{tau_bar}. To upper bound
the numerator $\P(|\rho_{n+1}|^\prime >h,|\rho_n|^\prime\le h)$ we
first write
\be\lbl{num_bound}\P(|\rho_{n+1}|^\prime >h,|\rho_n|^\prime\le h)\le
\P(|y_{n+1}|^\prime>h-h_1,|y_n+x_n|^\prime\le h)+\P(|x_{n+1}|^\prime>h_1).\ee 
Now, using  that $|y_n|^\prime-|x_n|^\prime\le |y_n+x_n|^\prime$ we see that $|y_n+x_n|^\prime\le
h$ implies that either $|x_n|^\prime>h_1$ or  $|y_n|^\prime\le
h+h_1$.  Therefore, 
\begin{eqnarray*}&&\P(|y_{n+1}|^\prime>h-h_1,|y_n+x_n|^\prime\le h)\le
\P(|y_{n+1}|^\prime>h_0,|y_n|^\prime\le
h+h_1)+\P(|x_{n}|^\prime>h_1)\\&&
\quad  \le\P(|y_{n+1}|^\prime>h_0,|y_n|^\prime\le
h_0)+\P(h-h_1<|y_n|^\prime\le h+h_1)+\P(|x_{n}|^\prime>h_1)
\end{eqnarray*}

Combining this with \eqref{num_bound} we obtain
\[\frac{\P(|\rho_{n+1}|^\prime>h,|\bar
 y_n|^\prime\le h)}{\P(|\rho_n|^\prime\le h)}\le
\frac{\P(|y_{n+1}|^\prime>h_0,|y_n|^\prime\le
 h_0)+\P(h-h_1<|y_n|^\prime\le
 h+h_1)+2\P(|x_n|^\prime>h_1)}{\P(|y_n|^\prime\le h_0)-\P(|x_n|^\prime>h_1)}.\]
Using the inequality $\frac s{t-u}\le \frac{s+u}t$ valid for all
$t>0$, $0\le s\le t$ and $0\le u<t$ 
we can bound the last expression by
\[
\frac{\P(|y_{n+1}|^\prime>h_0,|y_n|^\prime\le
h_0)+\P(h-h_1<|y_n|^\prime\le
h+h_1)+3\P(|x_n|^\prime>h_1)}{\P(|y_n|^\prime\le h_0)}.\]
Furthermore, since $y_n$ is a continuous random variable, $x\to
\P(|y_n|^\prime\le x)$ is uniformly continuous and since
$h_1=O(\sigma^\alpha h^\beta)$ and $\P(|y_n|^\prime\le h_0)=O(1)$ and $\P(|x_n|^\prime>h_1)=O(\sigma^{2(1-\alpha)}h^{2(2-\beta)})$ we can bound the last quantity by
\be\lbl{last_bound}
\frac{\P(|y_{n+1}|^\prime>h_0,|y_n|^\prime\le
h_0)}{\P(|y_n|^\prime\le h_0)}+O(\sigma^\alpha h^{\beta})+O(\sigma^{2(1-\alpha)}h^{2(2-\beta)})
.
\ee
To optimize the last two terms on the right hand side of (\ref{last_bound}), we
choose $\alpha=2/3$ and $\beta=4/3$, which makes them $O((\sigma h^2)^{2/3})$.
The first term on the right--hand side od (\ref{last_bound}), by Theorem~\ref{asymptotics},
is of much smaller order. This concludes the proof of Theorem~\ref{full}.

\section{Applications}\lbl{example}
\setcounter{equation}{0}

In this section, we illustrate the analysis in the previous sections
with applications to three problems in mathematical neuroscience. 
To this end, we use two conductance-based models of single neurons
and a model of electrically coupled network of pancreatic $\beta-$cells.
In the absence of noise, all these models exhibit stable periodic 
oscillations. Adding noise to these models results in more complex 
stochastic oscillatory regimes: irregular bursting and mixed-mode 
oscillations. We show that despite different mathematical formulations 
of these models and different forms of resultant oscillations,
emergent stochastic oscillatory patterns in these models are formed
due to the random perturbations of stable limit cycle oscillations.
We show numerically that the number of spikes in one burst in models
in \S\ref{ikm} and \S\ref{beta} and the number of small oscillations
in \S\ref{mmo} are distributed approximately geometrically
in accord with Theorem~\ref{full}.  

Numerical examples in \S\ref{ikm} and \S\ref{beta} appeared before 
in \cite{HM, MZ11a} respectively. We use them here to give the reader
a feeling for the range of possible applications of our analytical results.
The numerical example in \S\ref{mmo} is new. In somewhat different dynamical
regime for a stochastically forced $2D$ FitzHugh-Nagumo oscillator, 
asymptotically geometric distribution of small amplitude oscillations
was shown in \cite{BL} (see also \cite{MVE} for related results). 
Finally, we note that there is another mechanism for generating
irregular 
bursting 
and mixed-mode oscillations featuring geometric distribution. 
It does not involve 
random perturbations, but is based on certain properties 
of chaotic attractors in 
closely related neuronal models (see \cite{M09, MY}).

\begin{figure}
\begin{center}
{\bf a}\epsfig{figure=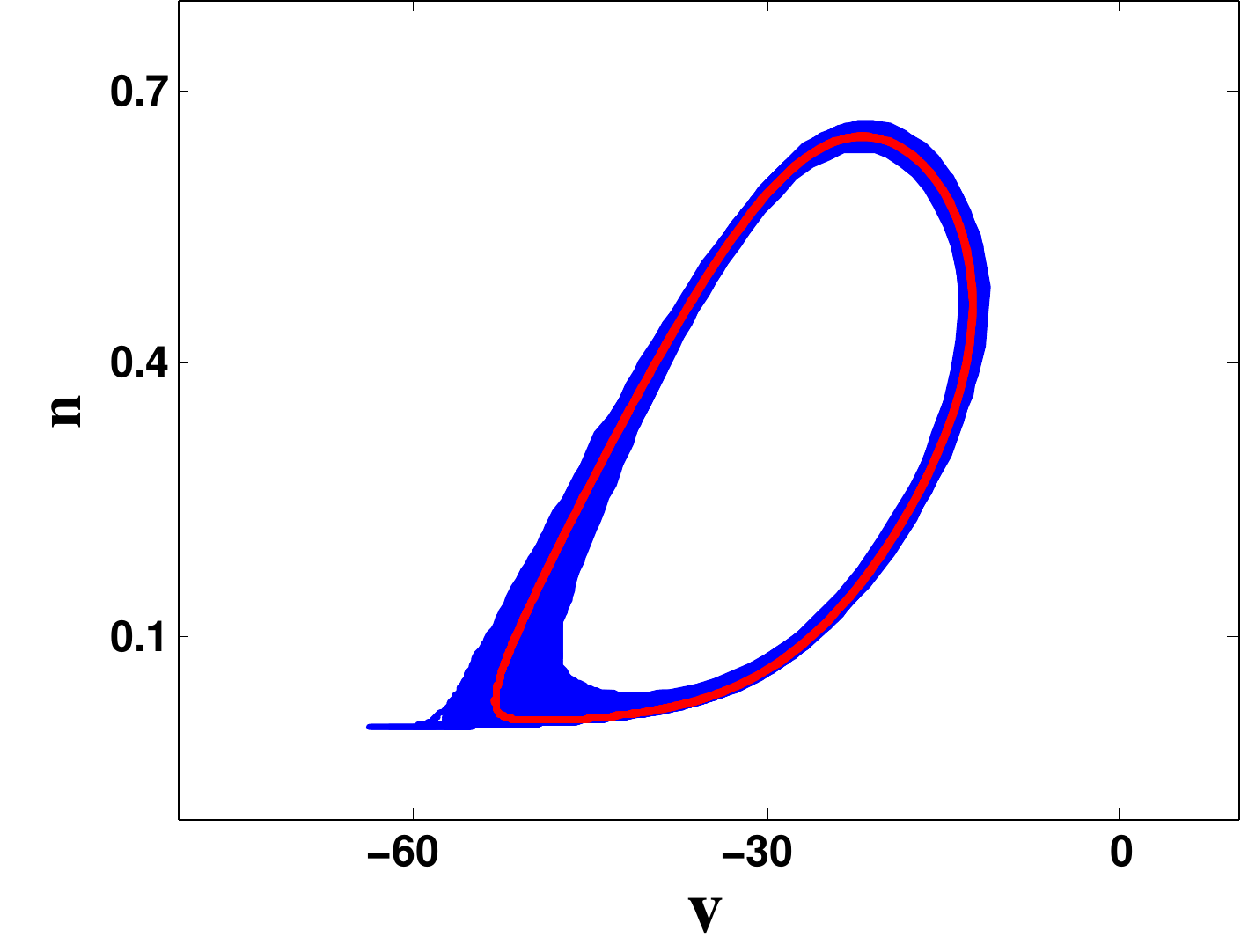, height=1.6in, width=2.05in}
{\bf b}\epsfig{figure=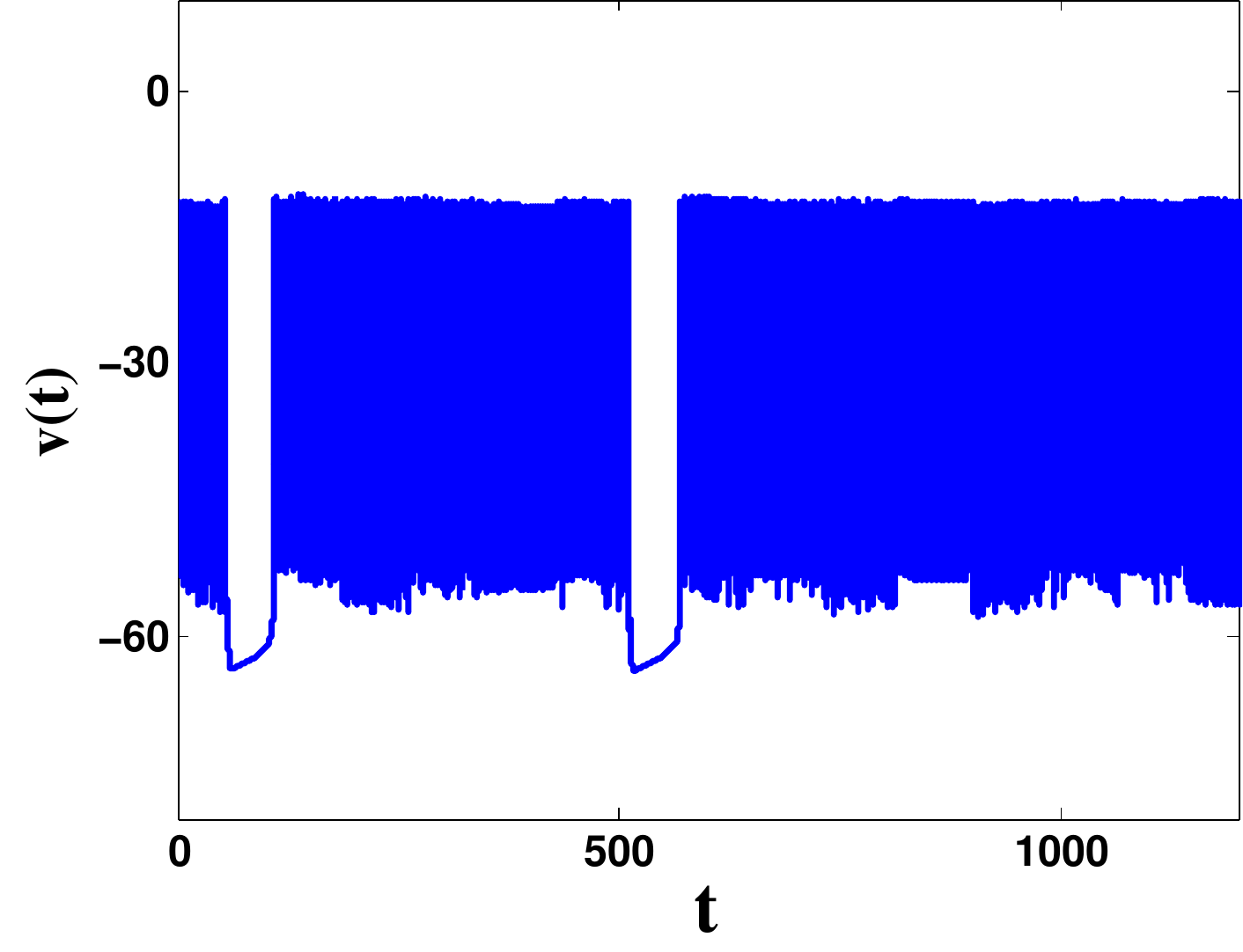, height=1.6in, width=2.05in}
{\bf c}\epsfig{figure=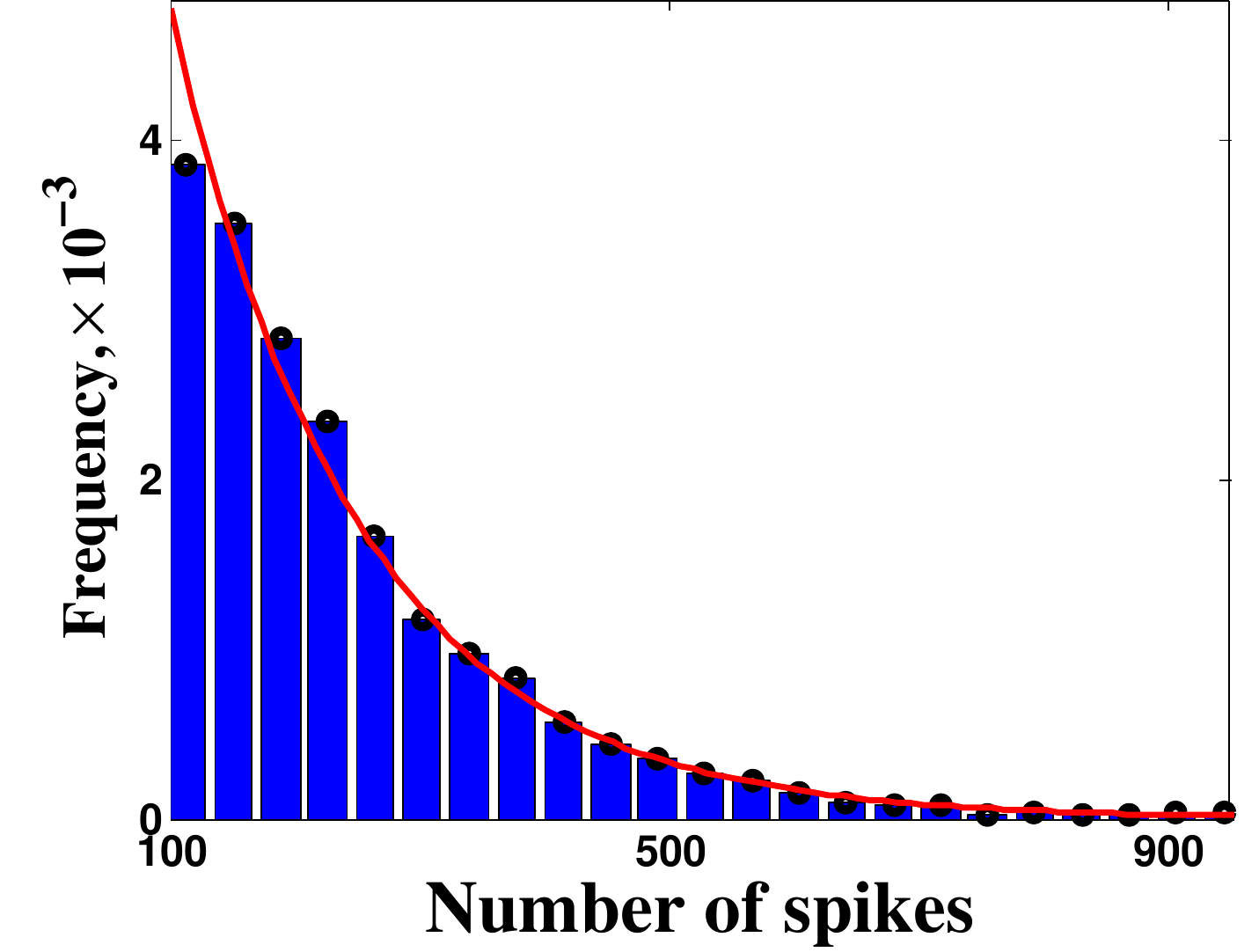, height=1.6in, width=2.05in}
\end{center}
\caption{Randomly perturbed model of neural oscillator
(\ref{10.1})-(\ref{10.3}). (a) The projection of a periodic
trajectory of the unperturbed deterministic system (red) and
that of the randomly perturbed one (blue) onto
$n$-$v$ plane. (b) The timeseries generated by the randomly 
perturbed model: series of fast oscillations closely resembling
oscillations of the unperturbed deterministic system are
separated by periods of no activity. (c) The tail of the normalized histogram
of the number of spikes in one burst is fitted with an exponential curve.
}
\lbl{f.1}
\end{figure}

\subsection{Bursting in a single cell model} \lbl{ikm}
Our first example concerns bursting. This is a common 
pattern of electrical activity in neurons. 
Bursting is characterized by series of fast oscillations separated
by periods of no activity (Fig.~\ref{f.1}b). Understanding dynamical
mechanisms underlying bursting is an important problem in
mathematical neuroscience \cite{Izh}. Below, we describe a 
mechanism of stochastic bursting based on random perturbations
of a stable limit cycle.

The following system of three differential equations describes the dynamics
of  a hypothetical neuron (cf.\S 9.1.1, \cite{Izh}):
\begin{eqnarray}\lbl{10.1}
C_m\dot v &=&-g_{NaP}m_\infty(v)(v-E_{NaP})-g_{K} n(v-E_K)-g_{KM} y(v-E_{K})\\
\nonumber
&-&g_L(v-E_L)+I+ \sigma\dot w_t,\\
\lbl{10.2}
\dot n &=& {n_\infty(v)-n\over\tau_n},\\
\lbl{10.3}
\dot y &=& {y_\infty(v)-y\over\tau_y}.
\end{eqnarray}
Here, $v$, $n$, and $y$ stand for the cell membrane potential and two gating
variables respectively. Parameters $g_s$ and $E_s$ are the maximal conductance
and the reversal potential of different ionic currents $I_s$, 
$s \in \left\{NaP, K, KM, L\right\}$ and $I$ is the applied current.
The time constants $\tau_n$ and $\tau_y$ determine the rates of
activation in the populations of $K$ and $KM$ channels. The
steady-state functions are defined by
$
s_\infty (v) =\left(1+\exp\left(\frac{a_s-v}{b_s}\right)\right)^{-1},\; s\in\left\{m,n,y\right\}.
$
In addition, a white noise process of intensity $\sigma$ 
is added to the right hand side of the voltage equation to model various
fluctuations affecting membrane potential.

The parameter values are summarized in the following table.
\begin{center}
\textbf{Table}
\end{center}
\begin{center}
\begin{tabular}{|l|c||l|c||l|c||l|c||l|c|}
\hline
$g_{Na}$    &   20 $mS/cm^2$   & $g_K$    & 10 $mS/cm^2$ &  
$g_{KM}$    &   5 $mS/cm^2$ &   $g_L$ &  8 $mS/cm^2$ &  $E_{Na}$    & 60 $mV$ \\
$E_K$    &   -90 $mV$  & $E_l$    &   -80 $mV$ & 
$\tau_n$    & 0.152 $ms^{-1}$  &  $\tau_y$  &   20 $ms^{-1}$ & 
$\sigma$    &  1 \\
I & 5 $pA$ &  $a_m$ & -20 $mV$ & $a_n$ & -25 $mV$ &  $a_y$ & -10 $mV$
& $b_m$  & 15 \\
$b_n$ &5 &  $b_y$ &5 &  $C_m$& 1 $\mu F/ cm^2$ & & &\\
\hline
\end{tabular}
\end{center}
\begin{figure}
\begin{center}
{\bf a}\epsfig{figure=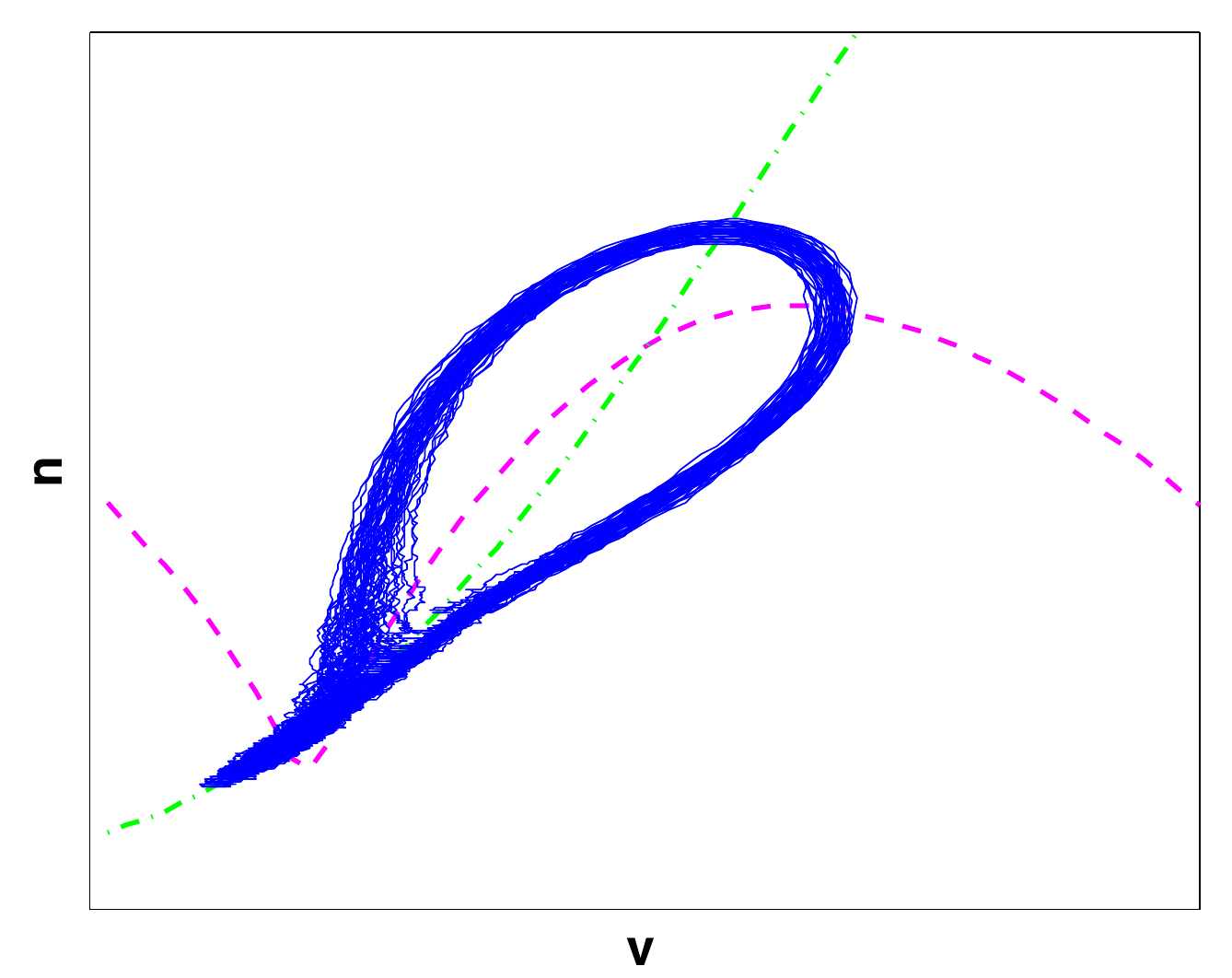, height=1.6in, width=2.05in}
{\bf b}\epsfig{figure=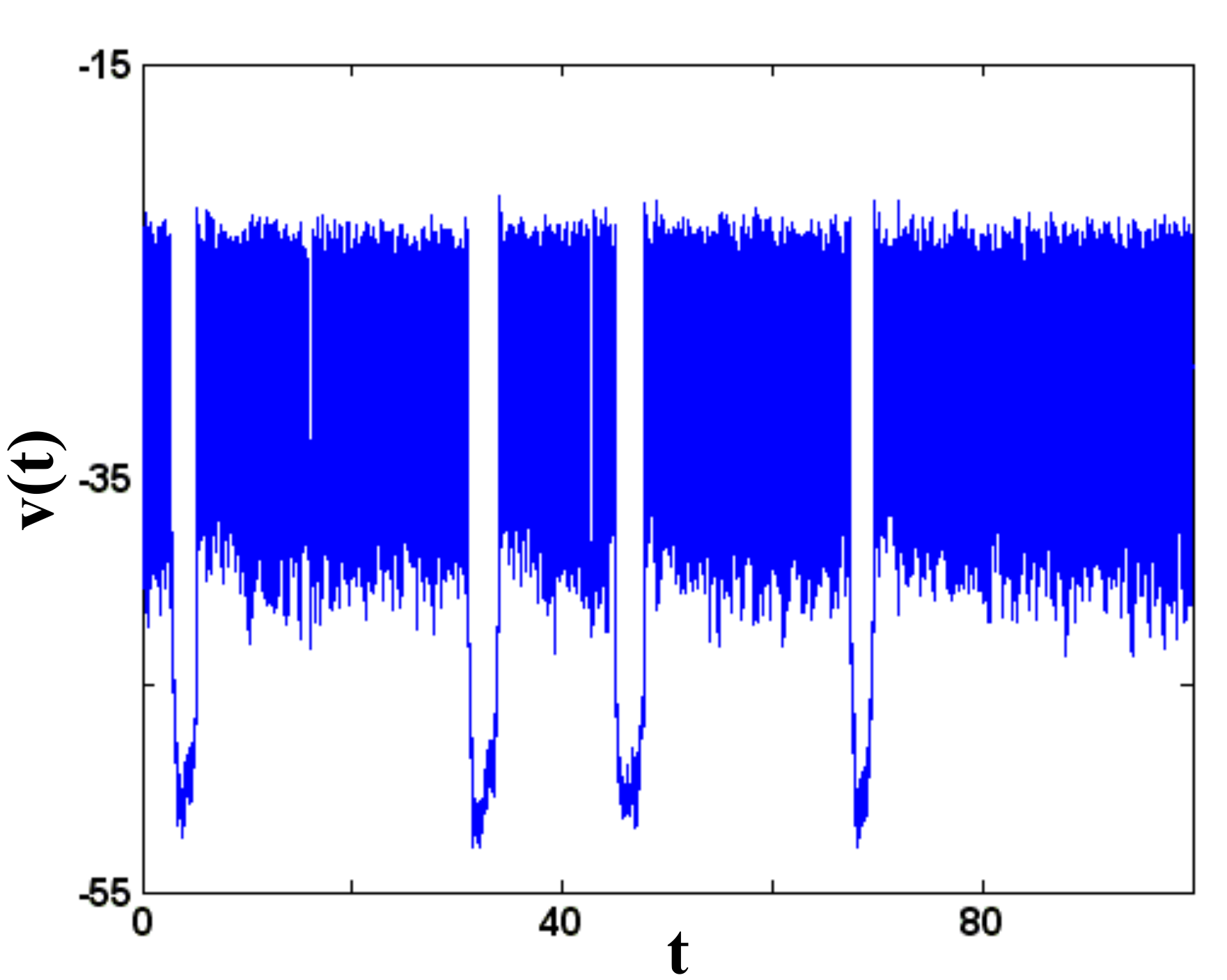, height=1.6in, width=2.05in}
{\bf c}\epsfig{figure=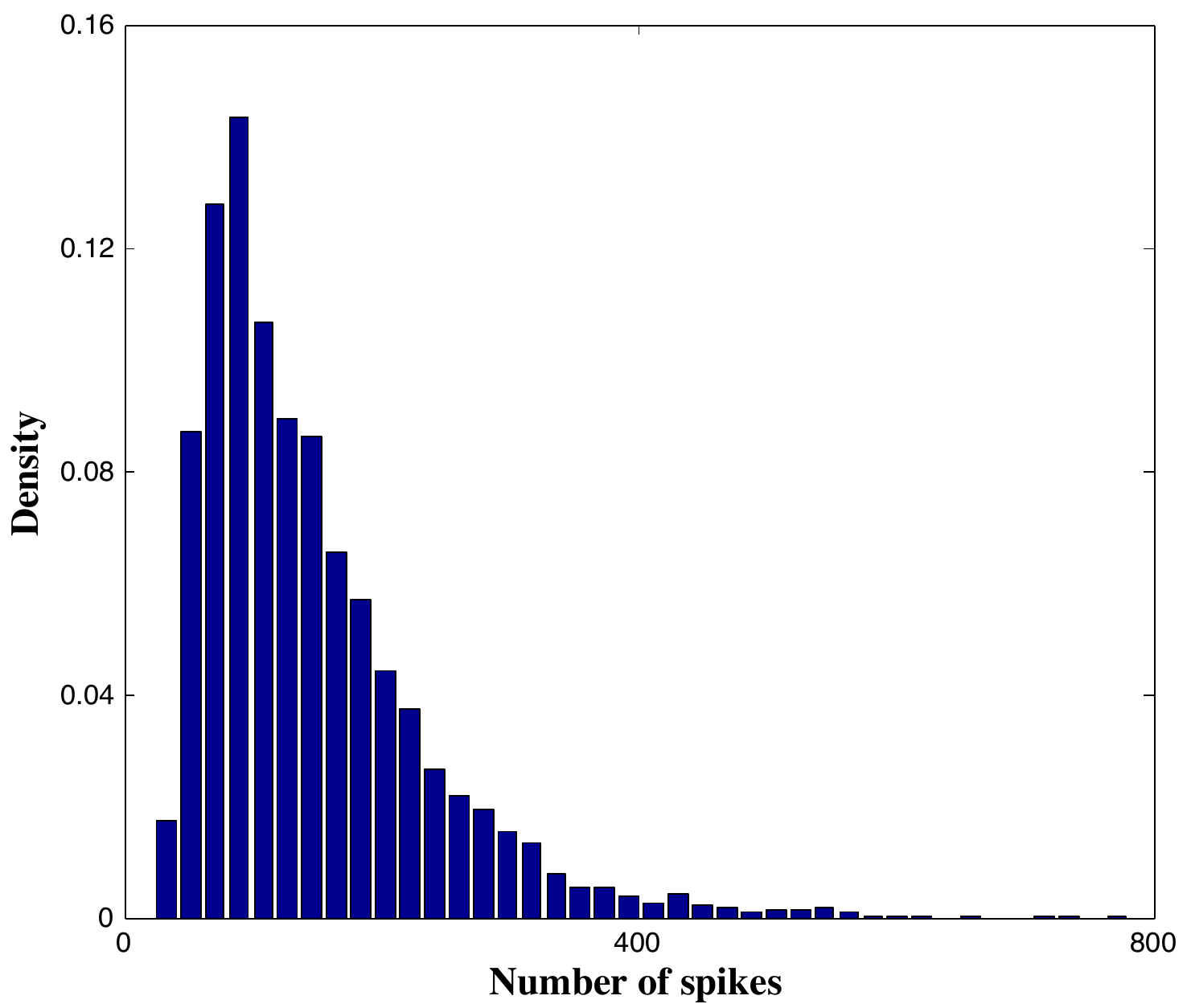, height=1.6in, width=2.05in}
\end{center}
\caption{
Randomly perturbed model of two coupled neural oscillators
(\ref{net.1})-(\ref{net.3}). (a) The projection of a periodic
trajectory of the unperturbed deterministic system (red) and that
of a randomly perturbed one (blue) onto
$n$-$v$ plane. (b) The timeseries generated by the randomly 
perturbed coupled model. (c) The normalized histogram
of the number of spikes in one synchronized burst 
features asymptotically geometric distribution.
}
\lbl{f.2}
\end{figure}

The deterministic model (\ref{10.1})-(\ref{10.3}) with $\sigma=0$
has a stable limit cycle $L$ (see Fig.~\ref{f.1}a). Under small random perturbation
($\sigma>0$), the trajectory of (\ref{10.1})-(\ref{10.3}) after a random number
of rotations around the limit cycle of the deterministic system leaves
the basin of attraction of the limit cycle and undergoes a large excursion 
in the phase space before coming back to a small neighborhood of $L$.
By applying the results of this paper to the  Poincare map of the limit cycle 
for this problem, one can show that the number of spikes in one burst is distributed
approximately geometrically. The tail of the distribution of the number
of spikes in one burst shown in  Fig.~\ref{f.1}c is well fitted by an exponential
function. 

Stochastic bursting in conductance-based models similar to (\ref{10.1})-(\ref{10.3})
was analyzed in \cite{HM}. The results of the paper provide a more general
treatment of  noise-induced bursting. Furthermore, the examples in 
the remainder
of this section lie outside the scope of applicability of the method in \cite{HM},
but can be effectively analyzed using the results of this paper. 

\subsection{Bursting in a coupled network} \lbl{beta}

Our next example presents a model of a coupled network generating synchronous
bursting oscillations. This is a model of an ensemble of pancreatic $\beta-$
cells coupled electrically. For simplicity, we use a two-cell network and refer
the interested reader to \cite{MZ11a} for the description and numerical simulations
of bigger networks.

 The dynamics of Cell~$i$ is governed by the following system of differential equations
\begin{eqnarray}\lbl{net.1}
C_m\dot v^{(i)} &=& -I_{ion} (v^{(i)}, n^{(i)}, y^{(i)})+ 
g (v^{(j)}-v^{(i)}) 
+\sigma \dot w^{(i)},\\
\lbl{net.2}
\dot n^{(i)} &=& {n_\infty(v^{(i)})-n^{(i)}\over \tau(v^{(i)})},\\
\lbl{net.3}
\dot y^{(i)} &=& \epsilon (I_{Ca}(v^{(i)})-ky^{(i)}), \quad i\in\{1,2\}, j=3-i.
\end{eqnarray}
The biophysical meaning of these equations is similar to that of the model in \S\ref{ikm}:
$v^{(i)}$, $n^{(i)}$, and $y^{(i)}$ stand for the membrane potential, gating variable,
and calcium concentration corresponding to Cell~$i$. Random processes $w^{(i)}, i\in\{1,2\},$ 
are  independent copies of standard Brownian motion.

The first term on the right hand side of (\ref{net.1}) models 
the combined effect of sodium and calcium currents, $I_{Na+Ca}$,
the calcium-dependent potassium current, $I_{KCa}$, delayed
rectifier $I_K$ and a small leak current, $I_l$
\be\lbl{apA.1} 
I_{ion}= I_{Na+Ca}+I_{KCa}+I_K+I_l,
\ee
where
\begin{eqnarray*}
I_{Na+Ca} &=& g_{Na+Ca}m_{\infty}(v)^3h_{\infty}(v)(E_I-v),\\
I_{KC} &=& {g_{KCa}u \over 1+u}(E_K-v),\\
I_K &=& g_Kn^4(E_K-v),\\
I_l &=& g_l(E_l-v).
\end{eqnarray*}
The steady state functions used to model the ionic currents
above are given by   
$$
f_\infty(v)=\frac{\alpha_f (v)}{\alpha_f(v)+\beta_f(v)},\quad f\in \{m,h,n\},
$$
where
\begin{eqnarray*}
\alpha_m = 0.1(v+25)(1-\exp\{-0.1(v+25)\})^{-1},  \beta_m=4\exp\{-(v+50)/18\}, & 
\alpha_h = 0.07\exp\{-0.05(v+50)\},\\
\beta_h=(1+\exp\{-0.1(v+20)\})^{-1},
\alpha_n = 0.01(v+20)(1-\exp\{-0.1(v+20)\})^{-1}, & 
\beta_n=0.125\exp\{ -(v+30)/80\}.
\end{eqnarray*}
The time constant of the delayed rectifier is given by
$$
\tau_n = (230(\alpha_n+\beta_n))^{-1}.
$$
The values of the remaining parameters are summarized in the following table.
\begin{center}
\textbf{Table}
\end{center}
\begin{center}
\begin{tabular}{|l|c||l|c||l|c||l|c|}
\hline
$g_{Na+Ca}$    &   1800$s^{-1}$   & $E_{Na+Ca}$    & 100$mV$ & $g_K$    &   1700$s^{-1}$    & $E_K$    & -75$mV$ \\
$k_C$    &   $\frac{\{2,12\}}{18}$$mV$  & $g_l$    &   7$s^{-1}$    & $E_l$    & -40$mV$ & $g_{KC}$    &   12$s^{-1}$  \\
$E_{Ca}$    & 100$mV$ & $\epsilon$    &     0.03$mV^{-1}s^{-1}$ & $C_m$ & $1\mu F/cm^2$ &  $\sigma$ & 
10 $mVs^{-1}$\\
\hline
\end{tabular}
\end{center}

When uncoupled ($g=0$), the models of each cell generate stable oscillations similar
to those generated by the model in \S\ref{ikm}. One can show that for sufficiently
strong coupling ($g\gg 1$), the deterministic coupled system ($\sigma=0$) has 
a stable limit cycle corresponding to synchronous oscillations in both
cells  (cf.~\cite{M11}). In the presence of noise ($\sigma>0$), the trajectory of the coupled system
once in a while leaves the basin of attraction of the limit cycle. This terminates
one burst. Since one spike in a burst corespond to one iteration of the PM, we 
conclude that the number of spikes in one synchronized burst of the coupled system
must be distributed approximatelly geometrically as shown in Fig.~\ref{f.2}.

\subsection{Mixed-mode oscillations}\lbl{mmo}
Our final example deals with 
mixed-mode oscillations, another type of nonlinear 
oscillations that are important in neuroscience \cite{DK, BL, MY, MVE}.
To this end, we use a modification of the Hodgkin-Huxley model
of a neuron in the regime close to the
Andronov-Hopf bifurcation. It was introduced by Doi and Kumagai in \cite{DK}.
 The model consists of the differential equations
for the membrane potential $v$ and two gating variables $n$ and $h$:
\begin{eqnarray}\lbl{Veqn}
C_m\dot v &=& -g_{Na}m_\infty(v)^3 h (v-E_{Na}) -g_k n^4 (v-E_K) -g_l(v-E_l)+\sigma\dot w,\\
\lbl{Neqn}    
\dot n &=& {n_\infty(v)-n\over \tau_n(v)},\\
\lbl{Heqn}    
\dot h &=& {h_\infty(v)-h\over \tau_h(v)}.
\end{eqnarray}

The biophysical meaning of these equations is similar to that of
the model discussed in \S\ref{ikm}. The steady state functions 
$$
{f}_\infty(v)=\frac{\alpha_f (v)}{\alpha_f(v)+\beta_f(v)},\quad f\in \{m,h,n\},
$$
and the time constants 
$$
\tau_f(v)=\frac{1}{\alpha_f(v)+\beta_f(v)},\quad f\in \{h,n\},
$$
are defined using 
\begin{eqnarray*}
\alpha_m = 0.1(v-25)(1-\exp\{0.1(25-v)\})^{-1},  \beta_m=4\exp\{-v/18\}, & 
\alpha_h = 0.07\exp\{-0.05 v\},\\
\beta_h= (1+\exp\{0.1(30-v)\})^{-1},
\alpha_n = 0.01(v-10)(1-\exp\{0.1(10-v)\})^{-1}, & 
\beta_n=0.125\exp\{ -v/80\}.
\end{eqnarray*}
The values of the remaining parameters are summarized in the following table.
\begin{center}
\textbf{Table}
\end{center}
\begin{center}
\begin{tabular}{|l|c||l|c||l|c||l|c||l|c|}
\hline
$g_{Na}$    &   120 $mS/cm^2$   & $g_K$    & 36 $mS/cm^2$ &   $g_l$    &   0.3$mS/cm^2$ &   $E_{Na}$    & 115 $mV$ &  $E_K$    &   -12 $mV$  \\
$E_l$    &   10.5999$mV$ & $\bar\tau_n$    & 20  &  $\bar \tau_h$  &   1 & 
$\sigma$    &  $5\cdot 10^{-4}$ &  & \\
\hline
\end{tabular}
\end{center}

The parameters in the neuronal model (\ref{Veqn})-(\ref{Heqn})
are chosen such that the deterministic model ($\sigma=0$)
has a stable limit cycle, whose projection onto $v-n$ plane
is shown in Fig.~\ref{f.3}a. Under the action of noise,
after a random number of rotations around the periodic orbit of the
deterministic system, the trajectory leaves the vicinity of 
the limit cycle. The global structure of the vector field of the
deterministic model guarantees that the trajectory of the randomly
perturbed system returns to the neighborhood of the limit cycle
after each excursion in the phase space. This results in 
mixed-mode oscillations of the membrane potential consisting
of a random  number of small oscillations separated by large 
spikes (see Fig.~\ref{f.3}b). In accord with the results of this
paper, we find that the distribution of the number of small
oscillations generated by this model is approximately geometric
(see Fig.~\ref{f.3}c). 
\begin{figure}
\begin{center}
{\bf a}\epsfig{figure=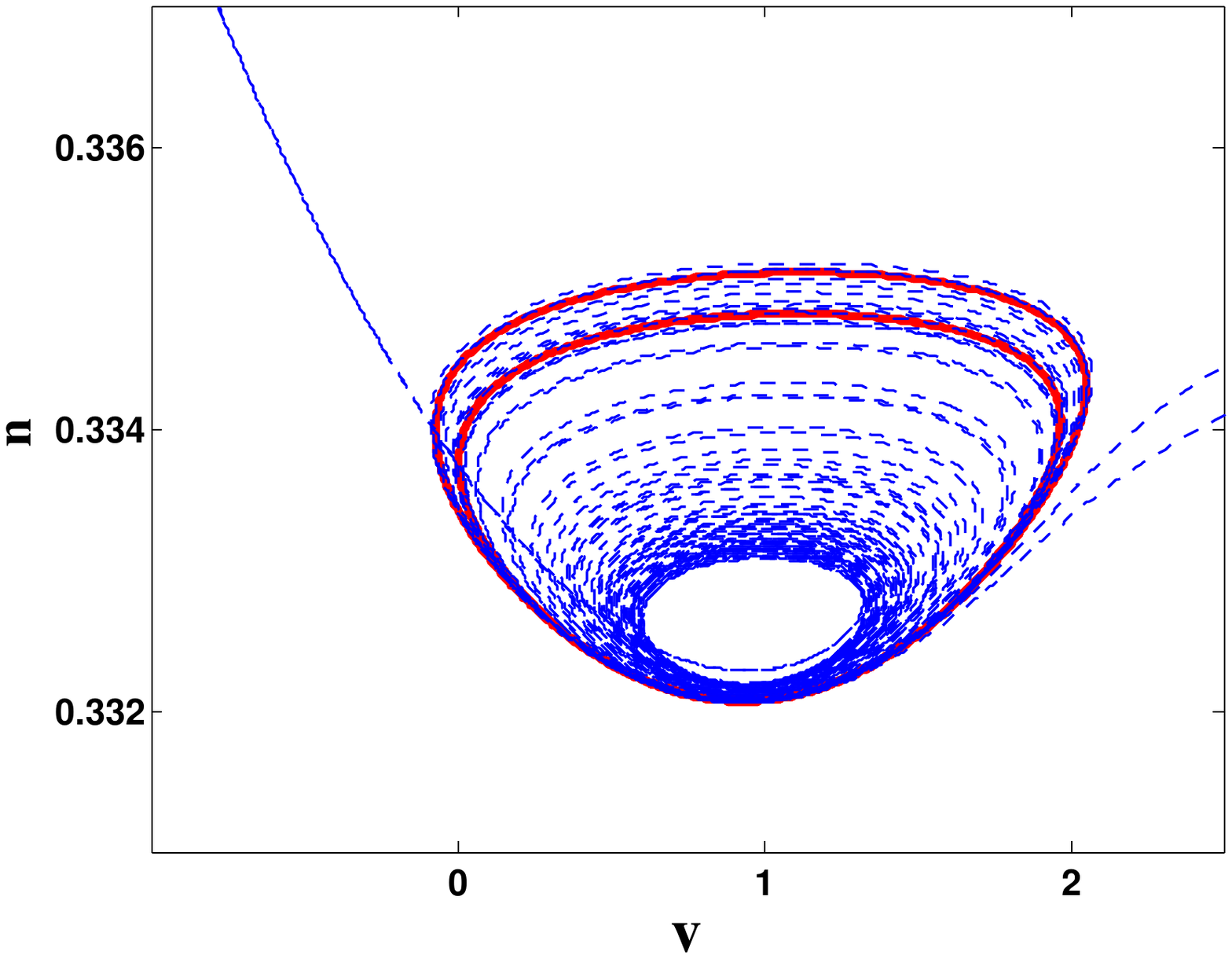, height=1.6in, width=2.05in}
{\bf b}\epsfig{figure=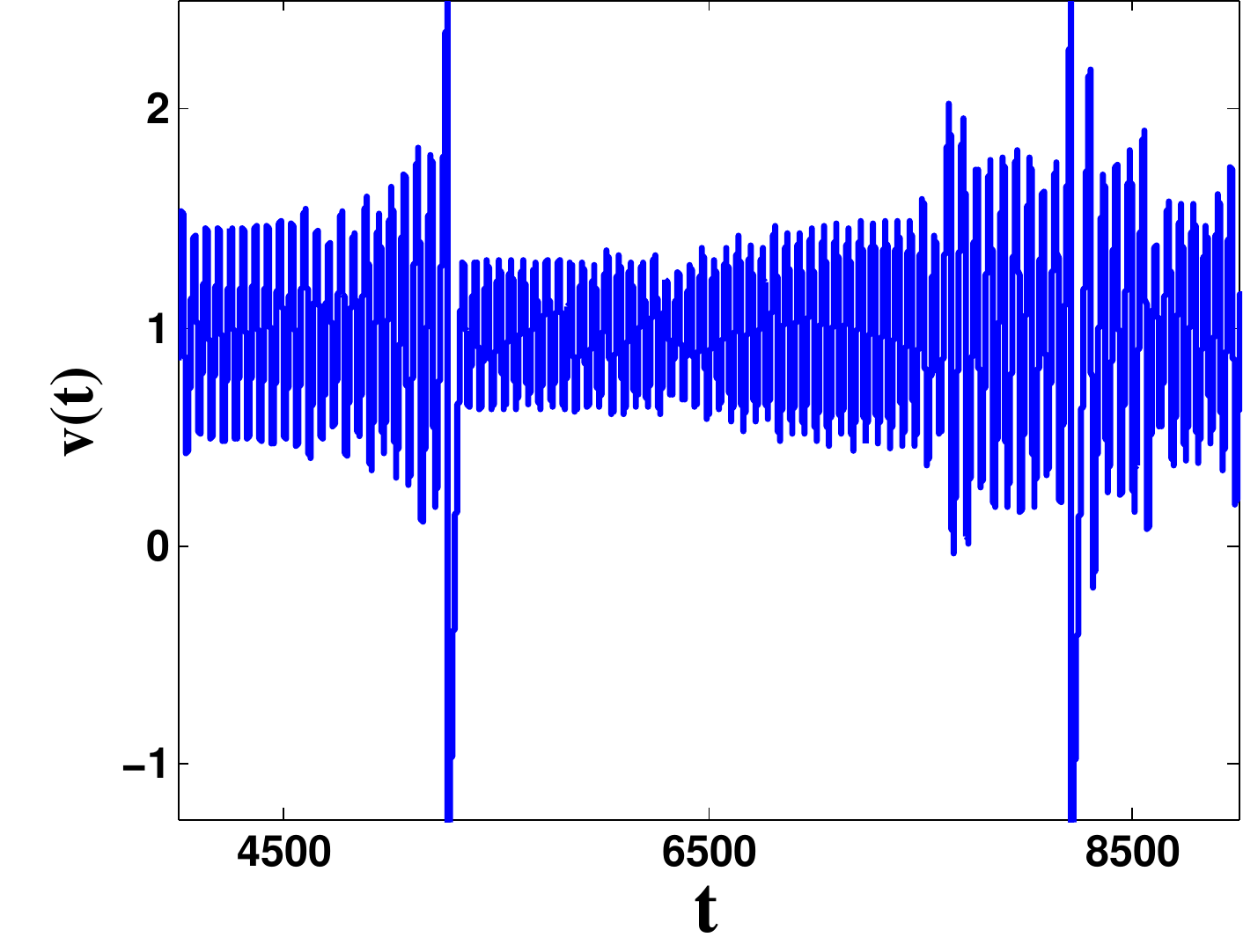, height=1.6in, width=2.05in}
{\bf c}\epsfig{figure=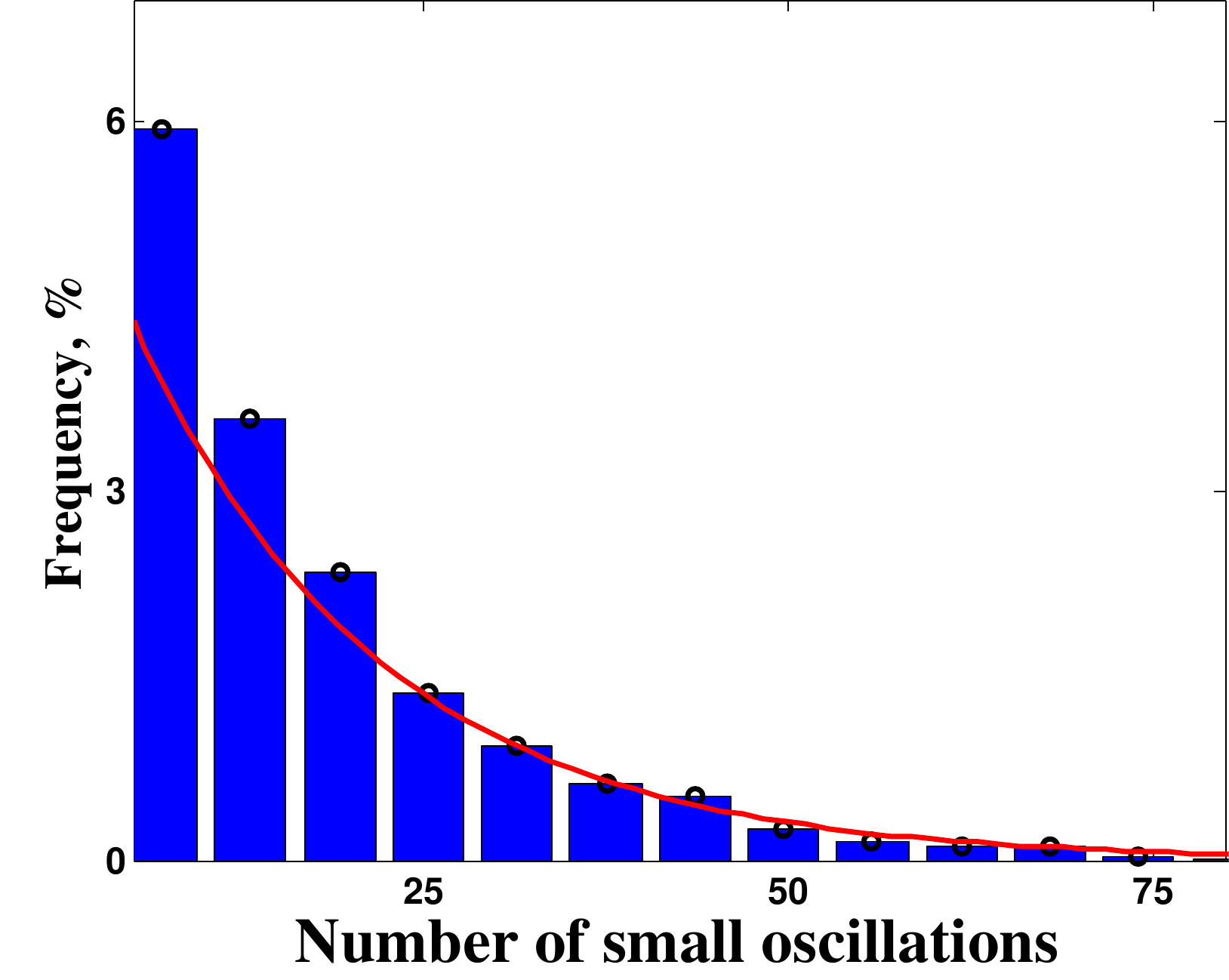, height=1.6in, width=2.05in}
\end{center}
\caption{
Noise-induced mixed-mode oscillations generated by (\ref{Veqn})-(\ref{Heqn}). 
(a) The projection of a periodic
trajectory of the unperturbed deterministic system (red) and that
of a randomly perturbed one (blue) onto
$n$-$v$ plane. (b) The timeseries generated by the randomly 
perturbed coupled model. (c) The normalized histogram
of the number of small oscillations between two consecutive spikes 
features asymptotically geometric distribution.
}
\lbl{f.3}
\end{figure}

\section{Discussion}\lbl{discuss}
\setcounter{equation}{0}
Since the time of Lyapunov and Poincare, problems in
nonlinear oscillations stimulated and guided the development of 
the geometric theory of ordinary differential equations
and its diverse applications to physics and biology
\cite{AVK, hale_osc, Malkin-a}.  Many effective analytical
methods have been developed for studying periodic motion and
effects of forcing in deterministic systems \cite{VIA, CH82, hale_odes}.
In contrast, apart from the large deviation type techniques
\cite{FW} and those for slow-fast systems \cite{BG}, there are
few general analytical approaches available for studying effects 
of random forcing on nonlinear oscillations and the analysis
of such systems is often done on the case-by-case basis
\cite{AIN, ARR, BL, DNR, HM}.

In the geometric theory of differential equations, 
the principal tool for studying
stability of periodic motion is the reduction to a PM.
This work presents a systematic construction of the  PM for an important class
of randomly perturbed problems: a limit cycle oscillator
forced by small white noise. For  trajectories of 
so-obtained PM,
we analyzed the statistics of the first exit times from a 
small neighborhood of the origin, which corresponds to the
periodic solution of the unperturbed system. We showed
that if the periodic solution of the deterministic system
is asymptotically stable the first exit times have approximately 
geometric distribution. This result implies universality
of the geometric distribution in diverse oscillatory regimes
generated by random perturbations of stable limit cycles.

Irregular oscillations featuring this dynamical mechanism
are common for differential equation models in applied science and,
in particular, in mathematical biology.
We provided three representative examples from biophysics: irregular 
bursting and mixed-mode oscillations generated by randomly 
perturbed conductance-based models of neurons and synchronous
noise-induced bursting in randomly forced neuronal network.
We showed that these
dynamical regimes are caused by random perturbations of stable limit cycles
and, therefore, all of them feature geometric distribution.

\noindent {\bf Acknowledgements.} 
This work was partially supported by a grant  from the Simons
Foundation  (grant no. 208766 to PH) and an NSF grant 
(DMS 1109367 to GM).

\end{document}